\documentclass[a4paper, 10pt]{amsart}
\pdfoutput=1

\usepackage{amssymb, amsmath}
\usepackage[english]{babel}
\usepackage{ifthen}
\usepackage{ifpdf}
\usepackage{graphicx}
\usepackage{url}
\usepackage{amscd}
\usepackage{a4wide}
\usepackage{todonotes}

\usepackage[pdftex, plainpages = false, 
    pdfpagelabels,  pdfpagelayout = SinglePage, 
    bookmarks = true, bookmarksopen = true, 
    bookmarksnumbered = true, linktocpage,
    colorlinks = true, 
    linkcolor = blue, urlcolor  = blue, 
    citecolor = red,  anchorcolor = green, 
    hyperindex = true, hyperfigures]{hyperref} 

\DeclareGraphicsExtensions{.png, .jpg, .pdf}
\graphicspath{{img/}{img/}{img/}}

\newcommand{\ol}{\overline}
\def\ds{\displaystyle}

\newcommand{\bpr}{\begin{proof}}
\newcommand{\epr}{\end{proof}}
\newcommand{\rank}{\operatorname{rank}}

\theoremstyle{plain}
\newtheorem{theorem}{Theorem}[section]
\newtheorem{lemma}[theorem]{Lemma}
\newtheorem{corollary}[theorem]{Corollary}
\newtheorem{proposition}[theorem]{Proposition}

\theoremstyle{definition}
\newtheorem{definition}[theorem]{Definition}
\newtheorem{example}[theorem]{Example}
\newtheorem{exercise}[theorem]{Exercise}
\newtheorem{conjecture}[theorem]{Conjecture}

\theoremstyle{remark}
\newtheorem{remark}[theorem]{Remark}

\numberwithin{equation}{section}

\newcommand{\bt}[1]{\begin{theorem}\label{#1}}
\newcommand{\bc}[1]{\begin{corollary}\label{#1}}
\newcommand{\bl}[1]{\begin{lemma}\label{#1}}
\newcommand{\bp}[1]{\begin{proposition}\label{#1}}
\newcommand{\be}[1]{\begin{example}\label{#1}}
\newcommand{\bd}[1]{\begin{definition}\label{#1}}
\newcommand{\br}[1]{\begin{remark}\label{#1}}
\newcommand{\bx}[1]{\begin{exercise}\label{#1}}
\newcommand{\bcon}[1]{\begin{conjecture}\label{#1}}

\newcommand{\et}{\end{theorem}}
\newcommand{\ec}{\end{corollary}}
\newcommand{\el}{\end{lemma}}
\newcommand{\ep}{\end{proposition}}
\newcommand{\ee}{\end{example}}
\newcommand{\ed}{\end{definition}}
\newcommand{\exc}{\end{exercise}}
\newcommand{\er}{\end{remark}}
\newcommand{\econ}{\end{conjecture}}

\def\A  {\mathcal{A}}
\def\FA {\mathcal{F(A)}}
\def\Ch {\mathcal{C(\A)}}

\def \F {\mathcal{F}}
\def \B {\mathcal{B}}
\def \P {\mathcal{P}}

\def \G {\mathcal{G}}
\def \s {\mathcal{S}}
\def \L {\mathcal{L}}

\def \h{\mathcal{H}}

\def\R{\mathbb{R}}
\def\N{\mathbb{N}}
\def\C{\mathbb{C}}
\def\Z{\mathbb{Z}}
\def\p{\mathbb{P}}

\begin{document}

\title{Arrangements of spheres and projective spaces}
\author{Priyavrat Deshpande}
\address{Chennai Mathematical Institute\\ India}
\email{pdeshpande@cmi.ac.in}
\date{}

\begin{abstract}
We develop the theory of arrangements of spheres. Consider a finite collection of codimension-$1$ subspheres in a positive-dimensional sphere. 
There are two posets associated with this collection: the poset of faces and the poset of intersections. 
We also associate a topological space: the complement of the union of tangent bundles of these subspheres in the tangent bundle of the ambient sphere. 
We call this space the tangent bundle complement. 
As in the case of hyperplane arrangements the aim of this new notion is to understand the interaction between the combinatorics of the intersections and the topology of the tangent bundle complement. 
In the present paper we find a closed form formula for the homotopy type of the complement and express some of its topological invariants in terms of the associated combinatorial information. 
\end{abstract}

\maketitle

\section{Introduction}\label{ch1}

An arrangement of hyperplanes is a finite set $\A$ consisting of codimension-$1$ subspaces of $\R^l$. 
These hyperplanes and their intersections induce a polyhedral stratification of $\R^l$. 
The combinatorial information of an arrangement $\A$ is contained in two posets, namely, the \textit{face poset} which consists of all the strata and the \textit{intersection poset} which contains all possible intersections of hyperplanes in $\A$. 
A topological space associated with $\A$, denoted $M(\A)$, is the complement of the union of the complexified hyperplanes in $\C^l$.
It is an open submanifold of $\C^l$ with the homotopy type of a finite-dimensional CW complex \cite[Section 5.1]{orlik92}. 
The study of this complement was initiated in the works of Fox and Neuwirth, Arnol'd, Brieskorn and Deligne in the $60$'s and $70$'s (see \cite[Section 5.1]{orlik92}). 
One of the aspects of the theory of arrangements is to understand the interaction between the combinatorial data of an arrangement and the topology of $M(\A)$. 
For example, the cohomology ring of the complement, known as the \textit{Orlik-Solomon algebra} is completely determined by the intersection data \cite[Section 5.4]{orlik92}. 
A pioneering result by Salvetti in \cite{sal1} states that 
the homotopy type of the complement is determined by the face poset.  \par

A generalization of hyperplane arrangements was introduced by the author in \cite{deshpande_thesis11} where a study of arrangements of codimension-$1$ submanifolds in a smooth manifold was initiated. In this paper we focus on a particular example: arrangements of spheres. 
Given a smooth sphere $S^l$ we consider a finite collection of codimension-$1$ sub-spheres, denoted by  $\A$, which satisfy reasonably nice conditions.
For example, these sub-spheres are tamely embedded, their intersections are hyperplane-like and they induce a stratification of the ambient sphere such that all the strata are contractible. 
Consequently, one can define face and intersection posets in this context.  
The topological space associated with such a collection is the complement of the union of tangent bundles of these sub-spheres in $TS^l$. 
We call this space the \textit{tangent bundle complement} and denote it by $M(\A)$.
We ask the same question, to what extent does the combinatorics of $\A$ help determine the topology of $M(\A)$ ?\par

We explore this interaction of combinatorics and topology by first describing a regular cell complex that has the homotopy type of $M(\A)$.
The construction of this complex relies on the order relations in the face poset and is a generalization of the classical Salvetti complex.
We then concentrate only on those arrangements which exhibit certain antipodal symmetry.
For these so-called \textit{mirrored arrangements} we find a closed form formula for the homotopy type of $M(\A)$.
We then show that the cohomology groups of $M(\A)$ are determined by the intersection data. 
Moreover, the coholomogy ring of $M(\A)$ can be expressed as a direct sum of an Orlik-Solomon algebra and a free abelian group in the top dimension. The rank of this top-dimensional free abelian group is equal to the number of graded pieces in the Orlik-Solomon algebra.
We also identify a class of arrangements for which the word problem for $\pi_1(M(\A))$ is solvable.
\par
In case of mirrored arrangements, as a consequence of the antipodal symmetry, we can define \textit{projective arrangements}, i.e., a finite collection of subspaces homeomorphic to $\R \p^{l-1}$ in $\R \p^l$. 
We exploit this antipodal symmetry to its full extent and derive similar results regarding the tangent bundle complement.
For example, the antipodal map helps understand the homotopy type of the tangent bundle complement as well as its fundamental group. 

\par 
An important motivation to study hyperplane arrangements comes from their natural connection with the Coxeter groups and the associated Artin groups.  
Let $W$ be a finite, irreducible Coxeter group of rank $n$. 
It acts linearly (in fact as origin-fixing isometries) on a real vector space $V$ of dimension $n$. Such a group is generated by reflections and has the following presentation - 
\[W = \langle s_1, \dots, s_n~|~ s_i^2 = 1, (s_is_j)^{m_{ij}} = 1,\quad \forall i\neq j \hbox{~and~} 2\leq m_{ij} < \infty\rangle. \] 
Its action on $V$ is not free; each reflection in $W$ fixes a hyperplane. 
The union of these reflecting hyerplanes is the \textit{reflection arrangement}, denoted $\A_W$, associated to $W$. 
The complement of these fixed hyperplanes is a disjoint union of open simplicial cones called (Weyl) chambers. 
Under the $W$ action these chambers are permuted freely (see \cite[Chapter 6]{davisbook08}). \par 

Complexifying this situation we get a finite arrangement of complex hyperplanes in $V\otimes\C$.
The complement of the union of these hyperplanes, denoted $M_W$, is connected and admits a fixed point free action of $W$. 
Brieskorn \cite{bries73} showed that the fundamental group of the orbit space $N_W$ has the following presentation -
\[\langle s_1,\dots, s_n\mid \underbrace{s_is_js_i\cdots}_{m_{ij}} = \underbrace{s_js_is_j\cdots}_{m_{ij}} \quad \forall i\neq j\rangle. \]
This group is known as the Artin group associated to $W$ and is denoted by $A_W$. 
There is a natural surjection from $A_W$ onto $W$ whose kernel is the so-called pure Artin group $PA_W$. 
It is the fundamental group of $M_W$. 
If $W$ is the symmetric group (i.e., type $A$ Coxeter group) then $A_W$ is the braid group and $PA_W$ is the pure braid group.\par  

Deligne \cite{deli72} showed that the universal cover of $N_W$ is contractible. 
Hence $N_W$ is a $K(A_W, 1)$ space. 
Subsequent study of these groups is much influenced by Deligne's work. Some of the important properties of Artin groups were proved by expanding on his ideas, notably the biautomatic nature of these groups \cite{charney92}. Simply put, it says that the Artin groups have solvable word and conjugacy problem. We refer the reader to \cite[Section 1.2]{charney07} for details and references.\par

We investigate sphere arrangements with a similar motivation. 
It is well known that the finite subgroups of isometries of a sphere generated by reflections are in fact Coxeter groups (see \cite[Chapter 10]{davisbook08} and \cite{gutkin86}).
Each reflection in this Coxeter transformation group fixes a codimension-$1$ subsphere giving rise to a sphere arrangement. The complement of this arrangement is a disjoint union of `spherical' simplices and they are freely permuted by the action. 
Since the group acts via isometries the action extends to the tangent bundle of the sphere. The complement of the union of the tangent bundles of the fixed sub-spheres serves as the analogue of the space $M_W$ introduced above. The Coxeter transformation group acts fixed point freely on this complement. The fundamental group of the orbit space is the desired generalization of Artin groups. The main aim of this paper is to lay topological foundations for the study of these ``Artin-like" groups. 
We illustrate with an example.

\be{egfinal}
Consider the $1$-sphere $S^1$. In this case a Coxeter transformation group $W$ is a dihedral group of order $2n$ with the presentation $\langle r, s ~|~ r^2 = s^2 = (rs)^n = 1\rangle$. The $n$ reflections in $W$ fix $n$ $0$-spheres i.e., $2n$ points. Declaring one of the chambers as the fundamental chamber all others can be labeled by elements of $W\setminus\{1\}$. The $2n$ points in this arrangement can be labeled by conjugates of the two standard parabolic subgroups $W_r, W_s$ of $W$. See \cite[Chapter 5]{davisbook08} for details regarding such labeling.\par 

The tangent bundle complement is an infinite cylinder with $2n$ punctures. 
The Salvetti complex (see Section \ref{tbc} for its construction) has $2n$ $0$-cells with labels $\langle g, g\rangle$ for every $g\in W$. There are $4n$ $1$-cells with labels of the form $\langle hW', h \rangle$ where $W'$ is one of the standard parabolic subgroups and $h\in W$. 
The reader can verify that the boundary of this $1$-cell is $\{ \langle h, h\rangle, \langle g, g\rangle \}$ such that $g^{-1}h \in W'$. 
The `labeled' Salvetti complex inherits the free $W$-action on the tangent bundle complement. 
The orbit complex consists of exactly one $0$-cell and two $1$-cells with both their end points joined at the $0$-cell. 
It has the homotopy type of wedge of two circles.\par 

The generalized pure Artin group in this case is $F_{2n+1}$ the free group on $2n+1$ generators. The generalized Artin group is $F_2$ and we get the following exact sequence 
\[1 \rightarrow F_{2n+1} \hookrightarrow F_2 \twoheadrightarrow W \rightarrow 1.\]
If we were to denote the generators of $F_2$ as $r$ and $s$ then $F_{2n+1} = \pi_1(M(\A))$ has the following presentation - 
\[\langle r^2, s^2, (rs)^n, rs^2r, rsr^2sr,\dots, \underbrace{rsr\cdots}_{n-1}\epsilon^2 \underbrace{srs\cdots}_{n-1}, sr^2s,\dots, \underbrace{srs\cdots}_{n-1}\epsilon^2 \underbrace{rsr\cdots}_{n-1} \rangle \]    
where $\epsilon$ is $r$ or $s$ depending on the parity of $n$. 
\ee 
In a joint work with Ronno Das \cite{dasdesh14} we extend this correspondence to higher-dimensional spheres (in fact, to smooth manifolds).
In particular, we show that the homotopy equivalence between the tangent bundle complement and the Salvetti complex is $W$-equivariant. 
Further we explain how the group theoretic data can be used to define this complex.
We also describe a presentation for the associated groups. 
Current work in progress, among other things, focuses on computation of the (twisted) cohomology of $M(\A)$ with group ring coefficients.
\par 

The paper is organized as follows. Section \ref{sec1} is about the preliminaries of hyperplane arrangements. In Section \ref{sec2} we introduce the new objects of study, arrangements of spheres and the tangent bundle complement. In Section \ref{sec3} we look at how the combinatorics of intersections determines the topology of the complement. We investigate the fundamental group in Section \ref{sec4}. In Section \ref{sec5} we look at arrangements of projective spaces. 

\section{Arrangements of Hyperplanes} \label{sec1}
Hyperplane arrangements arise naturally in geometric, algebraic and combinatorial instances. In this section we formally define hyperplane arrangements and the combinatorial data associated with it in the setting that is most relevant to our work. 

\bd{def1} A (real) \emph{arrangement of hyperplanes} is a collection $\A = \{H_1,\dots,H_k\}$ of finitely many hyperplanes in $\R^l$, $l\geq 1$. \ed

An arrangement is called \textit{central} if the intersection of all the hyperplanes in $\A$ is non-empty. However, we allow our arrangements to be non-central. For a subset $X$ of $\R^l$, the \emph{restriction} of  $\A$ to $X$ is the subarrangement $\A_X := \{H\in\A ~|~ X\subseteq H\}$.\par 

Analogously one can define hyperplane arrangements in $\C^l$ which are called \textit{complex arrangements}. To every real arrangement $\A$ there is an associated complex arrangement $\A_{\C}$; for every $H\in \A$ there is a hyperplane $H_{\C}\in \A_{\C}$ with the same defining equations as $H$. In this paper we focus on (complexified) real arrangements of hyperplanes.\par 
 
Associated with $\A$ there are two posets containing important information about the arrangement, namely, the face poset and the intersection poset. 

\bd{def2} The \emph{intersection poset} $L(\A)$ of $\A$ is the set of all intersections of hyperplanes, including $\R^l$ itself as the empty intersection, ordered by reverse inclusion. \ed

The intersection poset is a ranked poset with the rank of an element being the codimension of the corresponding intersection. 
The \textit{rank of an arrangement} $\A$ is defined to be the rank of its intersection poset. 
An arrangement is said to be \textit{essential} if its rank equals the dimension of the ambient space; without loss of generality we will from now assume this to be the case.
In general $L(\A)$ a (meet) semilattice; it is a lattice if and only if the arrangement is central.  

\bd{def3} The \emph{face poset} $\FA$ of $\A$ is the set of all faces ordered by topological inclusion: $F\leq G$ if and only if $F\subseteq\overline{G}$. \ed

Codimension-$0$ faces are called \emph{chambers}. The set of all chambers will be denoted by $\Ch$. A chamber is called bounded if it is a bounded subset of $\R^l$. Two chambers $C$ and $D$ are adjacent if they have a common face in their closure. 

The topological space associated with a real hyperplane arrangement $\A$ is its \emph{complexified complement} $M(\A)$ which is defined as follows: 

\bd{def4} \[M(\A) := \C^l \setminus \ds (\bigcup_{H\in\A} H_{\C}) \] 
where $H_{\C}$ is the hyperplane in $\C^l$ with the same defining equation as $H\in \A$.\ed

\subsection{The Salvetti Complex} \label{sec:TheSalvettiComplex}
In \cite{sal1} Salvetti constructed a regular CW-complex which has the homotopy type of the complexified complement. The construction uses the ordering in $\FA$. \par 
Let $\A$ be a hyperplane arrangement in $\R^l$. We construct a regular $l$-complex, called the Salvetti complex and denoted by $Sal(\A)$, by first describing its cells. The $k$-cells, for $0\leq k\leq l$, of $Sal(\A)$ are in one-to-one correspondence with the pairs $[F, C]$, where $F$ is a codimension-$k$ face of $\A$ and $C$ is a chamber whose closure contains $F$.\par 
Since $Sal(\A)$ is regular all the attaching maps are homeomorphisms. Hence it is enough to specify the boundary of each cell. A cell labelled $[F_1, C_1]$ is contained in the boundary of another cell labelled $[F_2, C_2]$ if and only if $F_1 \geq F_2$ in $\FA$ and $C_1, C_2$ are contained in the same chamber of $\A_{F_1}$. Now we state the seminal result of Salvetti.

\begin{theorem}[Salvetti \cite{sal1}]\label{thm0} Let $\A$ be an arrangement of real hyperplanes and $M(\A)$ be the complement of its complexification inside $\C^l$. Then there is an embedding of $Sal(\A)$ into $M(\A)$. Moreover there is a natural map in the other direction which is a deformation retraction. \end{theorem}

The above construction is generalized by Bj\"orner and Ziegler in \cite{bz92} where authors give a CW-complex with the homotopy type of the complement of a complex subspace arrangement.

\subsection{Cohomology of the Complement} \label{sec:CohomologyOfMA}

We begin by associating a combinatorially defined algebra, called the \emph{Orlik-Solomon algebra}, to a (complex) hyperplane arrangement. \par

Let $\A$ be a hyperplane arrangement. For every $1\leq p\leq n$ call a $p$-tuple $S = (H_1,\dots, H_p)$ of hyperplanes to be \textit{independent} if $\dim(H_1\cap\cdots\cap H_p) = l - p$ and call it \textit{dependent} if the intersection is nonempty and its codimension is strictly less than $p$. Geometrically, the independence implies that the hyperplanes of $S$ are in general position. \par 

Let $E_1$ be the free $\Z$-module generated by the elements $e_H$ for every $H\in \A$. Define $E(\A)$ to be the exterior algebra on $E_1$ and let $\partial$ denote the differential in $E(\A)$. For a $p$-tuple $S$ of hyperplanes we denote by $\bigcap S$ the intersection of elements in $S$ and by $e_S$ we mean $e_{H_1}\wedge\cdots\wedge e_{H_p}$. Let $I(\A)$ denote the ideal of $E$ generated by 
\[ \{ e_S~|~ \bigcap S = \emptyset\}\cup \{\partial e_S~|~ S\hbox{~is dependent}\}.\]

\bd{def15}The \emph{Orlik-Solomon algebra} of a (complex) arrangement $\A$ is the quotient algebra $E(\A)/I(\A)$ and denoted by $A(\A)$.  \ed 

The following important theorem shows how cohomology of $M(\A)$ depends on the intersection poset. It combines the work of Arnold, Brieskorn, Orlik and Solomon. For details and exact statements of their individual results see \cite[Chapter 3, Section 5.4]{orlik92}.

\bt{thm2sec1} Let $\A = \{H_1, \dots, H_n\}$ be a hyperplane arrangement in $\C^l$. For $H\in \A$ choose a linear form $\alpha_H\in (\C^l)^*$, such that $\ker(\alpha_H) = H$. Then the integral cohomology algebra of the complement is generated by the classes 
\[\omega_H := \ds\frac{1}{2\pi}\frac{d\alpha_H}{\alpha_H}. \]
The map $\gamma\colon A(\A)\to H^*(M(\A), \Z)$ defined by \[\gamma(e_H)\mapsto \omega_H\]
induces an isomorphism of graded $\Z$-algebras.  \et

This theorem asserts that a presentation of the cohomology algebra of $M(\A)$ can be constructed from the data that are encoded by the intersection poset. Let us mention one more theorem that explicitly states the role of the intersection poset in determining the cohomology of the complement. In particular the result states that there is a finer grading of cohomology groups indexed by the intersections and the rank of each cohomology group is determined by the M\"obius function of the intersection poset (see \cite[Proposition 3.75, Lemma 5.91]{orlik92}).\par 

\bt{thm112}Let $\A$ be a nonempty complex arrangement. For $X\in L(\A)$ let $M_X$ denote the complexified complement of the restricted arrangement $\A_X$. There are following isomorphisms for each $k\geq 0$
\[\ds \theta_k\colon \bigoplus_{\rank X = k}H^k(M_X)\to H^k(M) \]
induced by the inclusions $i_X\colon M\hookrightarrow M_X$. Moreover the rank of each cohomology group is determined by the following formula 
\[\rank H^k(M) = \sum_{X\in L_k} (-1)^{\rank(X)} \mu(\R^l, X) \]
where $\mu$ is the M\"obius function of $L(\A)$.
\et

\section{Arrangements of Spheres}\label{sec2}
We now introduce arrangements of codimension-$1$ subspheres in a sphere. First we isolate essential properties of a hyperplane arrangement:
\begin{enumerate}
\item[(1)] there are finitely many codimension $1$ subspaces each of which separates $\R^l$ into two components;
\item[(2)] there is a polyhedral stratification of the ambient space and the face poset of this stratification has the homotopy type of the ambient space. 
\end{enumerate}
\br{rem1} Recall that associated to every poset there is an abstract simplicial complex known as the \textit{order complex}. A $k$-simplex of the order complex corresponds to a $k$-chain of the poset. By homotopy type of the poset we mean the homotopy type of the geometric realization of the associated order complex. \er

In this section we first generalize the above properties in the context of spheres. 
Then we compare our definition with the topological representation of oriented matroids. 
Finally, we look at the combinatorics of the sphere arrangements. 

\subsection{Codimension-$1$ tame subspaces of Spheres}\label{sec2p1}
We start with a generalization of the property $(1)$ above. By an $l$-sphere $S^l$ we mean a smooth, closed $l$-manifold homeomorphic to the unit sphere in $\R^{l+1}$. The $0$-sphere $S^0$ consists of two points and we assume that the empty set is the sphere of dimension $-1$.\par 
If $S$ is an $(l-1)$-sphere embedded in $S^l$ ($l\geq 2$) as a closed subset then $S^l\setminus S$ has two connected components. Hence codimension-$1$ subspheres generalize hyperplanes in this respect. In general codimension-$1$ subspheres in a sphere could be very difficult to deal with. For example, consider the Alexander horned sphere. It is an embedding of $S^2$ inside $S^3$ such that one of the connected component of the complement is not even simply connected. In order to avoid such pathological instances we restrict ourselves to a nice class. A codimension-$1$ subsphere $S$ of $S^l$ is said to be \emph{tame} (or \textit{locally flat}) if for every $x\in S$ there is a neighbourhood $U_x$ of $x$ in $S^l$ such that $(U_x, U_x\cap S)\cong (\R^l, \R^{l-1})$. 
For such a subsphere the following statements are equivalent(see \cite[Theorem 1.8.2]{rushing73}):

\begin{enumerate}
\item there exists a homeomorphism $h$ of $S^l$ onto the standard unit sphere in $\R^{l+1}$ such that $h(S)$ is the equator cut out by the coordinate hyperplane $x_{l+1} = 0$; 
\item $S$ is homeomorphic to a piecewise-linearly embedded $(l-1)$-subsphere;
\item the closure of each connected component of $S^l\setminus S$ is homeomorphic to the $l$-ball.
\end{enumerate}

\begin{figure}[!ht]
\begin{center}
\includegraphics[scale=0.20,clip]{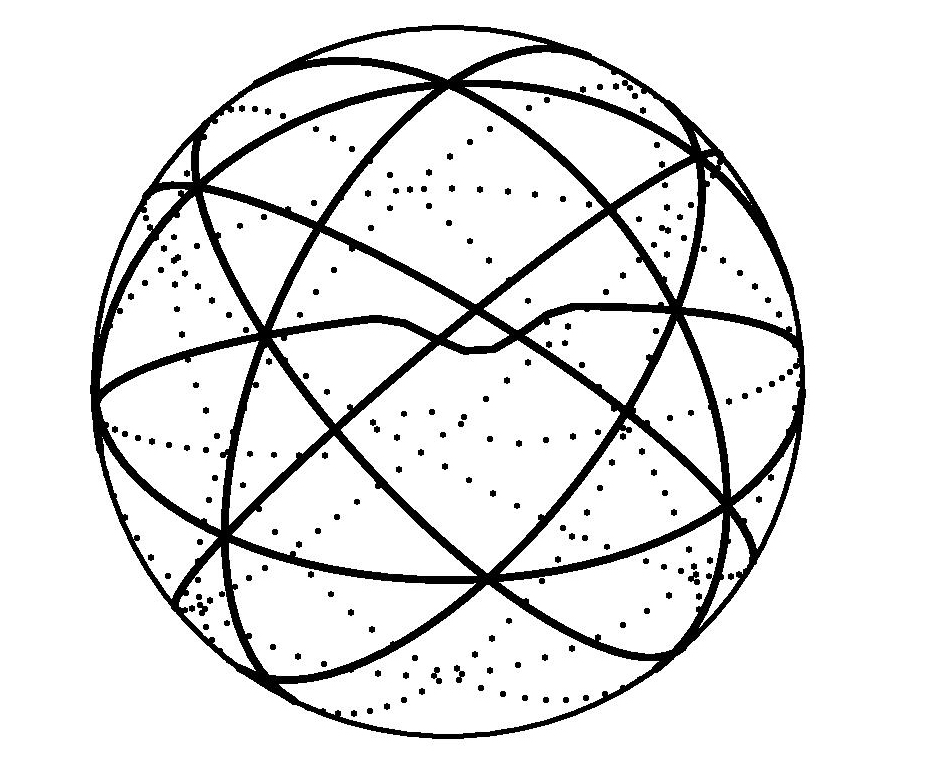}
\end{center}\caption{Non-Pappus arrangement}\label{nonpap}
\end{figure}

Tame subspheres need not intersect like hyperplanes. As an example, consider Figure \ref{nonpap}, it shows the non-Pappus arrangement on the unit sphere in $\R^3$. 
The cone over each of these $9$ circles is homeomorphic to a plane passing through the origin. Clearly there is no self-homeomorphism of $S^2$ such that the cone over the image of each of these circles is a $2$-dimensional subspace. This picture can arise as the boundary of a neighborhood of $2$-spheres intersecting in a $3$-sphere. We would like to avoid such situations as we are interested in dealing with the tangent bundle. \par 

We introduce a notion that will guarantee hyperplane-like intersections of subspheres. But first some notation. Let $\A= \{S_1, \dots, S_k \}$ be a collection of tame, codimension-$1$ subspheres in $S^l$. For every $x\in S^l$ and an open neighbourhood $V_x$ of $x$ homeomorphic to $\R^l$ let, 
\[ \A_x := \{S\cap V_x ~|~ x\in S\in \A \}.\] 
Denoting by $\bigcup \A_x$ we mean the union of elements of $\A_x$.

\bd{ch3def0} Let $\A = \{S_1,\dots,S_k\}$ be a collection of codimension-$1$, tame sub-spheres of $S^l$. We say that these sub-spheres have \emph{locally flat intersections} if for every $x\in S^l$ there exists an open neighbourhood $V_x$ and a homeomorphism $\phi\colon V_x\to \R^l$ such that $(V_x, \bigcup \A_x)\cong (\R^l, \bigcup \A')$ where $\A'$ is a central hyperplane arrangement in $\R^l$ with $\phi(x)$ as the common point. \ed

\subsection{Cellular stratification}
Now we generalize property (2). Let $\A = \{S_1,\dots,S_k \}$ be a collection of codimension-$1$, tame sub-spheres of $S^l$ with locally flat intersections. Let $\L$ be the set of all possible nonempty intersections of members of $\A$ and $\L^d$ be the subset containing codimension-$d$ intersections. We have $\bigcup \L^0 = S^l$ and $\bigcup \L^1 = \bigcup_{i=1}^k S_i$. For each $d\geq 0$ Consider the following subset of $S^l$ 
\begin{align*}
	\s^d(S^l) &=  \bigcup \L^d \setminus \bigcup \L^{d+1}.
\end{align*} 

Note that each $\s^i(S^l)$ may be disconnected and that the sphere can be expressed as the disjoint union of these connected components. We want these sets to define a `nice' stratification of $S^l$ hence we introduce the language of cellular stratified spaces developed by Tamaki in \cite{dai1}. 
Recall that a subset $A$ of a topological space $X$ is said to be \textit{locally closed} if every $x\in A$ has a neighborhood $U$ in $X$ with $A\cap U$ closed in $U$.

\bd{defn3s1}Let $X$ be a topological space and $\P$ be a poset. A \textit{stratification} of $X$ indexed by $\P$ is a surjective map $\sigma\colon X\to \P$ satisfying the following properties:
\begin{enumerate}
	\item For $p\in \mathrm{Im}\sigma$, $e_p := \sigma^{-1}(p)$ is connected and locally closed;
	\item for $p, q\in \mathrm{Im}\sigma,~ e_p\subseteq \ol{e_q} \iff p\leq q$;
	\item $e_p \bigcap \ol{e_q} \neq \emptyset \Longrightarrow e_p \subseteq \ol{e_q}$.
\end{enumerate}
The subspace $e_p$ is called the \textit{stratum with index $p$}.\ed

One can verify that the boundary of each stratum, $\partial e_p = \ol{e_p} - e_p$ is itself a union of strata. Such a stratification gives a decomposition of $X$. The indexing poset $\P$ is called the \emph{face poset} of the stratification. \par

It is now easy to check that the connected components of $\s^i(S^l)$ define a stratification of $X$. However it is not desirable to consider arbitrary stratifications. For example, in the case of two non-intersecting circles in $S^2$ there are three codimension-$0$ strata and two codimension-$1$ strata. But the resulting face poset does not have the homotopy type of the $2$-sphere. We need to focus on stratifications such that the strata are cells and the incidence relations between the strata recover the homotopy type of $S^l$. In order to achieve property (3) we assume that each stratum is a cell and the resulting stratification is a regular CW-complex.

\subsection{Definitions and examples}\label{sec2p3}
The desired generalization of hyperplane arrangements is the following: 
\bd{ch4def0}Let $S^l$ be a smooth sphere of dimension $l$. An \textbf{arrangement of spheres} is a finite collection 
$\A = \{S_1,\dots, S_k\}$ of codimension-$1$ smooth subspheres in $S^l$ such that:
	\begin{enumerate}
		\item the $S_i$'s have locally flat intersections (see Definition \ref{ch3def0});
		\item for all $I\subseteq \{1,\dots,k\}$ the intersection $\A_I := \bigcap_{i\in I} S_i$ is a sphere of some dimension;
		\item if $\A_I \nsubseteq S_j$, for some $I$ and some $j$, then $\A_I\cap S_j$ is a codimension-$1$ subsphere in $\A_I$; 
		\item the stratification induced by the intersections of $S_i$'s define the structure of a regular CW-complex.
	\end{enumerate} \ed

If there exists a fixed-point free, involutive diffeomorphism $f$ of the sphere such that for each $S\in \A$ we have $f(S) = S$ and $f(x)\neq x = f^2(x)~\forall x\in S^l$ then we call $\A$ a \emph{centrally symmetric arrangement of spheres}.

As in the case of hyperplane arrangements the combinatorial information associated with sphere arrangements is contained in the two posets which we now define.

\bd{def32}The \emph{intersection poset} denoted by $L(\A)$ is the set of connected components of all possible nonempty intersections of $S_i$'s ordered by reverse inclusion. By convention $S^l\in L(\A)$ as the least element.\ed

The intersection poset is a ranked poset. The rank of each element in $L(\A)$ is defined to be the codimension of the corresponding intersection.

\bd{def33}The intersections of these $S_i$'s in $\A$ define a stratification of $S^l$. The connected components in each stratum are called \emph{faces}. The collection of all the faces $\FA$ ordered by topological inclusion i.e., $F\leq G \Leftrightarrow F \subseteq \overline{G}$ is called the \emph{face poset}. The top-dimensional faces are called chambers and the set of all chambers is denoted by $\Ch$.\ed

For a face $F$ define its \textit{support} as the least-dimensional intersection containing $F$. The \textit{dimension of a face} is the dimension of its support. It is straightforward to see that the dimension function makes the face poset a ranked poset. \par 

\begin{figure}[!ht]
\begin{center}  \includegraphics[scale=0.5,clip]{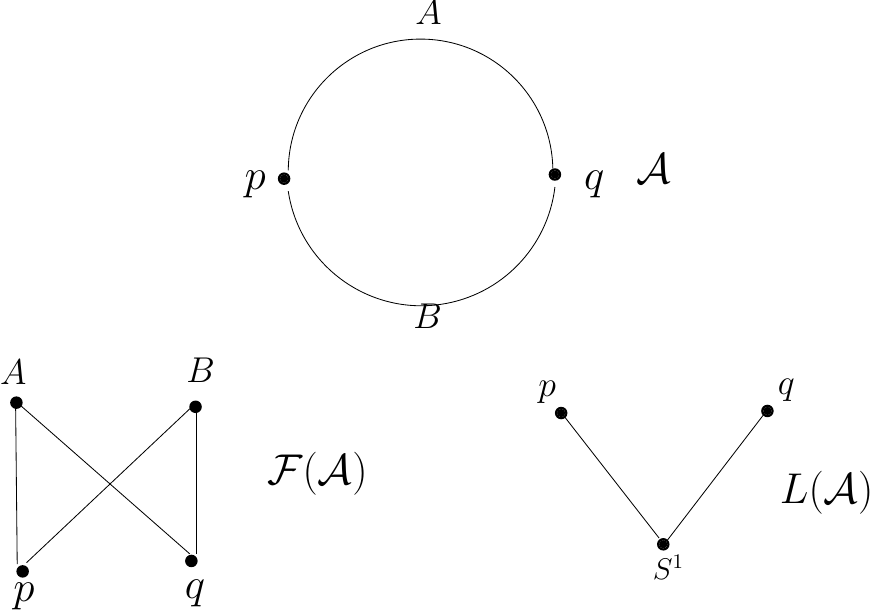} \end{center}
    \caption{Arrangement of $2$ points in a circle.}     \label{figdef1}     \end{figure}

\begin{figure}[!ht]
\begin{center}\includegraphics[scale=0.5,clip]{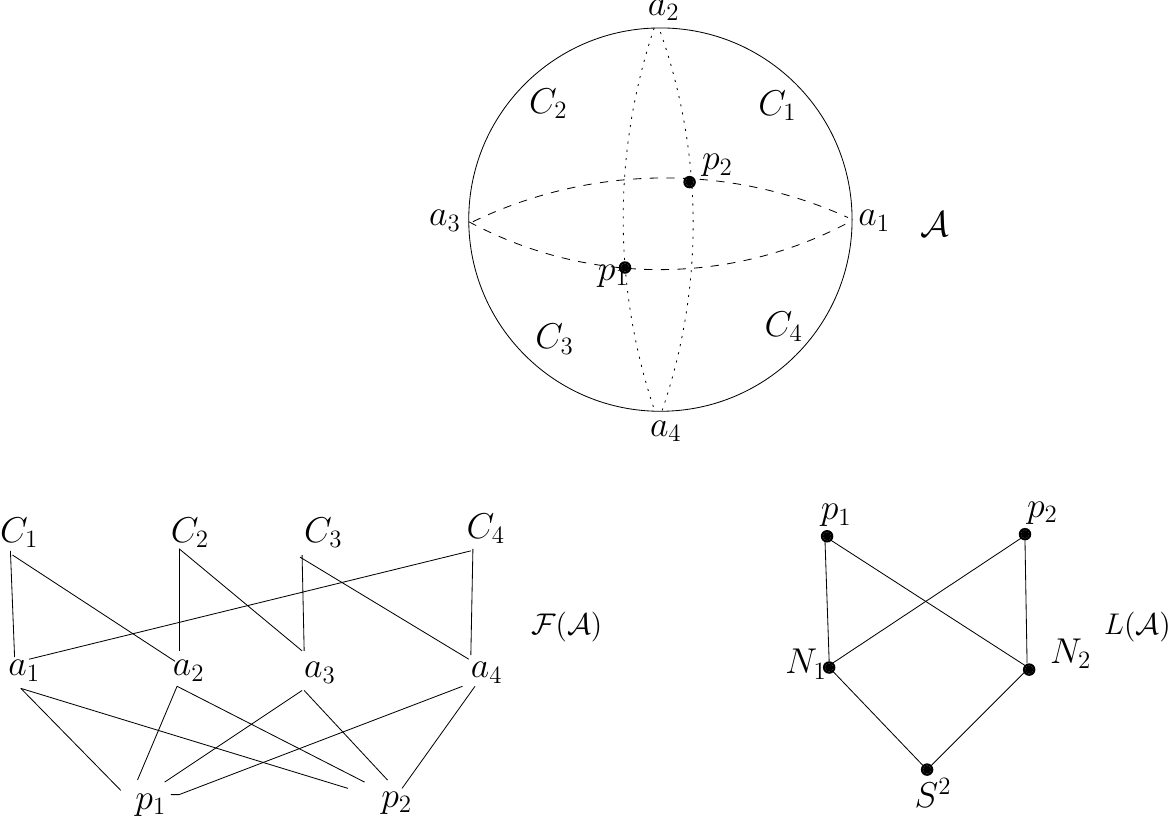}\end{center}
    \caption{Arrangement of $2$ circles in a sphere.} \label{sphere01}     \end{figure}

We now look at two examples of sphere arrangements. 

\be{ex41} Let $X$ be the circle $S^1$, a smooth one-dimensional manifold. The codimension-$1$ subspheres are the pairs of (diametrically opposite) points in $S^1$. Consider the arrangement $\A = \{p, q\}$ of $2$ such points. For both these points there is an open neighbourhood which is homeomorphic to an arrangement of a point in $\R$. Figure \ref{figdef1} shows this arrangement and the Hasse diagrams of the face poset and the intersection poset. 
\ee

\be{ex42}As a $2$-dimensional example consider an arrangement of $2$ great circles $N_1, N_2$ in $S^2$. Figure \ref{sphere01} shows this arrangement and the related posets. The face poset has two $0$-cells, four $1$-cells and four $2$-cells. Also note that the geometric realization of the face poset has the homotopy type of $S^2$. 
 \ee

\subsection{Topological representation of oriented matroids}
We now explore a connection between sphere arrangements and hyperplane arrangements using oriented matroids. The theory of oriented matroids is intimately connected with hyperplane arrangements. This combinatorial structure combines the information contained in face and intersection posets of a hyperplane arrangement. There are several (axiomatic) ways of defining oriented matroids. We refer the reader to the book of Bj\"orner et. al. \cite{ombook99} for various aspects related to oriented matroids. We do not intend to define and explain the properties of oriented matroids. Our aim is to compare their topological representation with the sphere arrangements.\par 

The oriented matroids which correspond to hyperplane arrangements are known as the \textit{realizable oriented matroids}.  There are oriented matroids that do not correspond to hyperplane arrangements (e.g., the non-Pappus configuration). Hence for a long time an important question in this field was to come up with the right topological model for oriented matroids. This was settled by Folkman and Lawrence in \cite{folk_lawr}. The Folkman-Lawrence Topological Representation Theorem states that in general oriented matroids correspond to certain collections of finitely many topological spheres and balls. These so-called pseudo-arrangements not only describe oriented matroids in the same way that $\R^l$ and collections of half spaces describe an obvious combinatorial structure but there is a one-to-one correspondence between such arrangements and the oriented matroids. \par 

In their original formulation Folkman and Lawrence introduced arrangements of pseudo-hemispheres.
Much simplification of their ideas was achieved by A. Mandel in his thesis \cite{mandel82}. He defined the notion of \textit{sphere systems} which we now state. 

\bd{defom2}A finite multi-set $\A = \{S_e\mid e\in E \}$ of codimension-$1$, tame subspheres in $S^l$ is called a \emph{sphere system} if the following conditions hold:
\begin{enumerate}
\item $S_A := \bigcap_{e\in A} S_e$ is a sphere, for all $A\subset E$.
\item If $S_A\nsubseteq S_e$ for $A\subset E, e\in E$, and $S^+_e$ and $S^-_e$ are the two sides of $S_e$, then $S_A\cap S_e$ is a subsphere in $S_A$ with sides $S_A\cap S^+_e$ and $S_A\cap S^-_e$.
\end{enumerate} \ed

A sphere system is said to be \textit{essential} if the intersection of all the sub-spheres is empty. It can be shown that the stratification of $S^l$ induced by an essential sphere system defines regular CW decomposition of the sphere \cite[Proposition 5.1.5]{ombook99}. The topological representation theorem states that (loop-free) oriented matroids of rank $l+1$ (upto reorientation and isomorphism) are in one-to-one correspondence with centrally symmetric, essential sphere systems in $S^l$. However, the sphere arrangements that we want to deal with are not general enough to represent oriented matroids. Note the differences between the definition of a sphere system and Definition \ref{ch4def0}. 
\begin{enumerate}
\item We assume that all intersections are locally flat,
\item the arrangement is \textit{repetition-free}, i.e., every subsphere appears exactly once.
\end{enumerate}

Given a central and essential arrangement of hyperplanes consider its intersection with the unit sphere; these intersections define a centrally symmetric sphere arrangement in the sense of Definition \ref{ch4def0}. However, the converse need not be true. Figure \ref{nonpap} shows an arrangement of $9$ pseudo-circles which is a sphere arrangement in the sense of Definition \ref{ch4def0} but it does not arise as an intersection with a central hyperplane arrangement. \par

The rank $3$ oriented matroids can be realized as centrally symmetric, repetition-free, essential sphere systems in $S^2$. The reader can verify that such sphere systems are arrangements of spheres (since intersections of pseudo-circles are locally-flat). 
Consequently, there is a one-to-one correspondence between  rank $3$ oriented matroids and the sphere arrangements in $S^2$. 
However, not all higher rank oriented matroids can be realized using sphere arrangements. For example, one can construct a sphere system in $S^3$ such that there exists at least one intersection whose spherical neighborhood looks like Figure \ref{nonpap}.

\br{remNew}We would like to clarify the distinction between Definition \ref{ch4def0} and Definition \ref{defom2}. 
The concept of sphere system is more general than that of an arrangement of spheres. 
It is clear that every sphere arrangement is a sphere system. 
However the converse is not true in general. 
There are examples of oriented matroids (say, of rank $4$) which do not correspond to any sphere arrangement in $S^3$.
In light of this observation it would be an interesting problem to obtain a combinatorial characterization of (non-realizable) oriented matroids which correspond to sphere arrangements.\er

A \textit{pseudo-hyperplane} is a tame embedding of a codimension-$1$ subspace in $\R^l$. 
Equivalently, it is the cone over a tame subsphere. An arrangement of pseudo-hyperplanes, intuitively, can be constructed by taking the cone over a sphere system. 

\bd{defpsh} A finite collection $\B = \{H_1,\dots, H_n\}$ of pseudo-hyperplanes in $\R^l$ is called an \textit{arrangement of pseudo-hyperplanes} if: \begin{enumerate}
\item For every $A\subseteq \{1,\dots, n\}$ the set $H_A := \bigcap_{i\in A}H_i$ is either empty or homeomorphic to some $\R^k$ for $0\leq k\leq l$.
\item For every $j\notin A$ either $H_A\subseteq H_j$, or $H_j\cap H_A$ is a locally flat embedding of a codimension-$1$ subspace of $H_A$.
\end{enumerate} \ed

We say that a pseudo-hyperplane arrangement is \textit{locally flat} if all the (nonempty) intersections are locally flat. 

The construction of the Salvetti complex and the Orlik-Solomon algebra hold true in case of pseudo-hyperplane arrangements. In fact, the Orlik-Solomon algebra associated to a pseudo-hyperplane arrangement is isomorphic to the cohomology algebra of the corresponding Salvetti complex (see \cite[Section 2.5]{ombook99} and \cite[Section 7]{bz92} for details). 


\subsection{Combinatorics of sphere arrangements}
We now take a closer look at the combinatorics of the incidence relations among the faces. Let $\A$ denote a sphere arrangement in $S^l$. A hypersphere $S$ in $\A$ is said to \emph{separate} two chambers $C$ and $D$ if they are contained in the distinct connected components of $S^l\setminus S$. For two chambers $C, D$ the set of all the hyperspheres that separate these two chambers is denoted by $R(C, D)$. The following lemma is now evident.

\bl{lem1s2c3}
Let $\A$ be an arrangement of spheres in $S^l$, an $l$-sphere. Let $C_1, C_2, C_3$ be three chambers of this arrangement. Then, 
\[R(C_1, C_3) = [R(C_1, C_2)\setminus R(C_2, C_3)] \cup [R(C_2, C_3)\setminus R(C_2, C_1)]. \] \el

The \emph{distance between two chambers} is defined as the cardinality of $R(C, D)$ and denoted by $d(C, D)$. Given a face $F$ and a chamber $C$ of a sphere arrangement $\A$ define the action of $F$ on $C$ as follows: 

\bd{def351} A face $F$ acts on a chamber $C$ to produce another chamber $F\circ C$ satisfying:
\begin{enumerate}
	\item $F \subseteq \overline{F\circ C}$,
	\item $d(C, F\circ C) = \hbox{~min~} \{d(C, C') ~|~ C'\in \Ch, F\subseteq\overline{C'} \}$. 
\end{enumerate} \ed

\bl{lem2s2c3}
With the same notation as above, the chamber $F\circ C$ always exists and is unique. \el

\bpr Clear. \epr

It is easy to check that if $C$ is a chamber and $F, F'$ are two faces such that $F' \geq F \hbox{~then ~} \ds F'\circ(F\circ C) = F'\circ C$. Moreover if $F\leq C$ then $F\circ C = C$.

\subsection{The tangent bundle complement}\label{tbc}
Recall that for a real hyperplane arrangement $\A$ the complexified complement $M(\A)$ is the complement of the union of complexified hyperplanes inside the complexified ambient vector space (Definition \ref{def4}). If one were to forget the complex structure on $\C^l$ then, topologically, it is just the tangent bundle of $\R^l$. Same is true for a hyperplane $H$ and its complexification $H_{\C}$. Hence the complexified complement of a hyperplane arrangement can also be considered as a complement inside the tangent bundle. We use this topological viewpoint to define a generalization of $M(\A)$ for sphere arrangements.

\bd{def35} As before let $\A =\{N_1,\dots, N_k \}$ be a sphere arrangement in $S^l$. Let $TS^l$ denote the tangent bundle of $S^l$ and let $T\A := \bigcup_{i=1}^k TN_i$. The \textbf{tangent bundle complement} of $\A$ is defined as \[M(\A) := TS^l\setminus T\A. \]  \ed

The above space was introduced in \cite[Chapter 3]{deshpande_thesis11} in the context of submanifold arrangements. We now construct a regular CW-complex, in the spirit of Salvetti's construction, that has the homotopy type of the tangent bundle complement. \par

We denote by $(S^l, \F(\A))$ the regular cell structure of $S^l$ induced by $\A$. 
We are interested in the dual cell structure which is obtained as follows.
For every face $F$ fix a point $x(F)\in F$ call it \textit{the barycenter} of $F$. 
Note that $\ol{F}$ is homeomorphic to an appropriate-dimensional disc $B_F$. 
Then there exists a regular cell structure of $B_F$ whose face poset is isomorphic to that of $\ol{F}$.	
For every $G < F$ the barycenter $x(G)$ determines a point $y_G$ of $B_F$.
Moreover, if $\gamma := G_0 < \cdots < G_k$ is a chain of faces of $F$ then form a simplex $\gamma_B$ of $B$ which is the convex hull of the vertices $y_{G_0}, \dots y_{G_k}$. 
Denote by $\Delta(\gamma)$ the image of $\gamma_B$ under the given homeomorphism. 
Note that $\Delta(\gamma)$ need not be the convex hull of $x(G_0), \dots x(G_k)$. 
Finally, denote by $F^*$ the union of all those $\Delta(\gamma)$'s which arise from chains ending in $F$ and call it the \textit{dual cell} of $F$.
The collection of all the dual cells defines a regular cell structure since link of each vertex is a sphere. 
We denote by $(S^l, \F^*(\A))$  this dual cell structure. 
Here $\F^*(\A)$ is the face poset of this cell structure with the partial order $\preceq$. 
Note that $\F^*(\A)$ is dual poset of $\F(A)$, i.e., $G^*\preceq F^* \iff F\leq G$.
Every $k$-face in $(S^l, \F(\A))$ corresponds to an $(l-k)$-cell in $(S^l, \F^*(\A))$ for $0\leq k\leq l$.  
\par

For the sake of notational simplicity we will denote the dual cell complex by $\F^*(\A)$ (and by $\F^*$ if the context is clear). Note that a $0$-cell $C^*$ is a vertex of a $k$-cell $F^*$ in $\F^*$ if and only if the closure $\overline{C}$ of the corresponding chamber contains the $(l-k)$-face $F$. The action of the faces on chambers that was introduced in Definition \ref{def351} is also valid for the dual cells. The symbol $F^*\circ C^*$ will denote the vertex of $F^*$ which is dual to the unique chamber closest to  $C$. \par

\noindent Given a sphere arrangement $\A$ in $S^l$ construct a regular $l$-complex $Sal(\A)$ as follows:\par

\noindent  The $0$-cells of $Sal(\A)$ correspond to $0$-cells of $\F^*$, which we denote by the pairs $\left<C^*; C^*\right>$.\\ 
For each $1$-cell $F^* \in \F^*$ with vertices $C^*_1, C^*_2$, take two homeomorphic copies of $F^*$ denoted by $\left<F^*; C^*_1\right>$ and $\left<F^*; C^*_2\right>$. Attach these two $1$-cells in $Sal(\A)_0$ (the $0$-skeleton) such that 
\[ \partial \left<F^*; C^*_i\right> = \{\left<C^*_1; C^*_1\right>, \left<C^*_2; C^*_2\right> \}\] 
for $i = 1, 2$. We put an orientation on the $1$-skeleton $Sal(\A)_1$ by directing each $1$-cell $\left<F^*; C^*\right>$ such that the initial vertex is $\left<C^*; C^*\right>$. \par

By induction assume that we have constructed the $(k-1)$-skeleton of $Sal(\A)$, $1\leq k-1 < l$. 
To each $k$-cell $G^* \in \F^*$ and to each of its vertex $C^*$ assign a $k$-cell $\left<G^*; C^*\right>$ whose face poset is isomorphic to that of $G^*$. 
Let $\phi(G^*, C^*)\colon \partial \left<G^*; C^*\right>\to Sal(\A)_{k-1}$ be the same characteristic map that identifies a $(k-1)$-cell $H^*\subseteq \partial G^*$ with the $k$-cell $\left<H^*; H^*\circ C^*\right>\subseteq \partial \left<G^*; C^*\right>$. Extend the map $\phi(G^*, C^*)$ to the whole of $\left<G^*; C^*\right>$ and use it as the attaching map, hence obtaining the $k$-skeleton. The boundary of every $k$-cell is given by 

\begin{equation}
\partial \left<F^*; C^*\right> = \bigcup_{G^* \prec F^*} \left<G^*; G^*\circ C^*\right>. \label{eq1s3c3}
\end{equation}

Now we state the theorem that justifies the construction of this cell complex. 

\bt{thm1c3s3} The regular CW-complex $Sal(\A)$ constructed above has the homotopy type of the tangent bundle complement $M(\A)$. \et

\bpr 
This is a special case of \cite[Theorem 3.3.7]{deshpande_thesis11}. 
We only sketch the proof here. 
The first step is to identify an open cover of the sphere indexed by the faces. There are three key properties that these open sets satisfy. First, for every face $F$ the corresponding open set $V_F$ is a regular neighborhood of $F$. Second, for another face $F'$ the intersection $V_F\cap F' \neq \emptyset$ if and only if $F\leq F'$. Finally, all these open sets and their non-empty intersections are contractible.
Now for any point on the manifold the tangent space at that point contains a arrangement of hyperplanes combinatorially equivalent to the local arrangement. Using the local trivialization one can construct an open covering of $M(\A)$ which is indexed by the pairs $\{(F, C)\in \F(\A)\times \Ch\mid F\leq C \}$ such that each of the open set is contractible and so are their intersections. This type of open covering satisfies the hypothesis of the nerve lemma. The final step is to establish the condition when two such open sets have a non-empty intersection. Thus providing an isomorphism between the nerve of this open cover and the face poset of the Salvetti complex constructed above. \epr


\be{ex1s3c3}As an example consider the arrangement of $2$ points in a circle (Example \ref{ex41}). The left side of the Figure \ref{circsal} below illustrates the arrangement with  the induced dual cell structure drawn using dotted lines. The right hand side shows the associated Salvetti complex with the cell labeling.  
\begin{figure}[!ht]
\begin{center}
	\includegraphics[scale=0.65,clip]{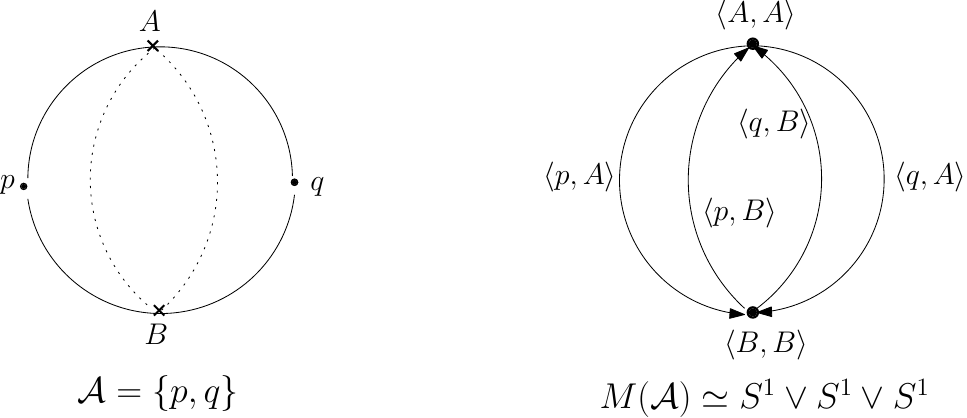}
\end{center}\caption{Arrangement in $S^1$ and the associated Salvetti complex}\label{circsal} 
\end{figure}\ee

We now look at some obvious properties of the above defined CW structure and also infer some more information about the tangent bundle complement.

\bt{thm2s3c3} Let $\A$ be a sphere arrangement in $S^l$ and let $Sal(\A)$ denote the associated Salvetti complex. Then
\begin{enumerate}
	\item there is a natural cellular map $\psi\colon Sal(\A)\to \F^*(\A)$ given by $\langle F^*, C^*\rangle\mapsto F^*$. The restriction of $\psi$ to the $0$-skeleton is a bijection and in general 
	\[ \psi^{-1}(F^*) = \{\langle F^*, C^*\rangle, C\in \Ch | C^* \preceq F^* \}.\] 
	\item  For every chamber $C$ there is a cellular map $\iota_C\colon \F^*(\A)\to Sal(\A)$ taking $F^*$ to $\langle F^*, F^*\circ C^*\rangle$ which is an embedding of $\F^*(\A)$ into $Sal(\A)$, and \[Sal(\A) = \bigcup_{C\in \Ch} \iota_C(\F^*).\] 
	\item The absolute value of the Euler characteristic of $M(\A)$ is the number of chambers.
\item Let $T\A$ denote the union of the tangent bundles of the submanifolds in $\A$ then,
\[\mathrm{rank}~\tilde{H}^i(TS^l, T\A) = \begin{cases} |\chi(M(\A))| &\hbox{~\emph{if}~} i = l \\ 0 &\hbox{~\emph{otherwise.}}\end{cases} \]
\end{enumerate} \et

\begin{proof}
Proofs of (1) and (2) are fairly straightforward. It follows that $S^l$ is homeomorphic to a retract of $M(\A)$. \par

We prove (3) by explicitly counting cells in the Salvetti complex. The Euler characteristic of a CW-complex is equal to the alternating sum of the number of cells of each dimension. Given a $k$-dimensional dual cell $F^*$ there are as many as $|\{C\in\Ch | F\leq C \}|$ $k$-dimensional cells in $Sal(\A)$. Hence for a $0$-cell $\langle C^*, C^*\rangle\in Sal(\A)$ the number of $k$-cells of $Sal(\A)$ with this particular vertex is equal to the number of $k$-cells of $\F^*(\A)$ that contain $C^*$. The alternating sum of number of cells that contain a particular vertex $C^*$ of $\F^*(\A)$ is equal to $1 - \chi(\mathrm{Lk}(C^*))$, where $\mathrm{Lk}(C^*)$ is the link of $C^*$ in $\F^*$. Applying this we get, 
\[\chi(Sal(\A)) = \sum_{C\in \Ch}(1 - \chi(\mathrm{Lk}(C^*))). \] 

Since $S^l$ is compact all the chambers are bounded we have $\mathrm{Lk}(C^*)\simeq S^{l-1}$. Thus, 
	\begin{align*}
				\chi(Sal(\A)) &= \sum_{C\in \Ch}(1 - \chi(\mathrm{Lk}(C^*)))  \\
						   &= \sum_{C \in \Ch}(1 - [1 + (-1)^{l-1}]) \\
						   &= (-1)^l \sum_{C \in \Ch}1.
	\end{align*}
Hence, \[ \chi(M(\A)) = (-1)^l \hbox{(number of chambers)}. \]
Let $\bigcup\A$ denote the union of hyperspheres in $\A$. Since $\A$ induces a regular cell decomposition $\bigcup\A$ has the homotopy type of wedge of $(l-1)$-spheres. The claim (4) follows from the homeomorphism of pairs $(TS^l, T\A) \cong (S^l, \bigcup\A )$.
\end{proof} 

\bc{lem1s6c3} Let $\A$ be a sphere arrangement in  $S^l, l\geq 2$. Then $M(\A)$ can not be an aspherical space. \ec

\br{remint} Let $Y$ be a positive-dimensional intersection of $\A$. The restriction of $\A$ to $Y$, i.e., the collection $\A_Y := \{S\cap Y\mid S\in\A, Y\nsubseteq S \}$ defines a sphere arrangement in $Y$. The reader can check using the inclusion $\iota_C$ (its restriction to $\F(\A_Y)$) from Theorem \ref{thm2s3c3} that $Y$ is homeomorphic to a retract of $Sal(\A)$.\er

\section{Topology of the Complement}\label{sec3}
The aim of this section is to investigate how the combinatorics of the associated posets affects the topology of the tangent bundle complement. 
Our investigation is based on a simple observation; if for a given centrally symmetric sphere arrangement there is an equator generically intersecting the sub-spheres then the restriction of the arrangement to the two hemispheres (i.e., components of the complement of the equator) gives combinatorially identical pseudo-hyperplane arrangements. 
We claim that these restricted pseudo-arrangements play a central role in understanding the topology of the complement. We identify a class of sphere arrangements for which it is easy to derive a closed form formula for the homotopy type of the complement. Then we establish a connection between the intersection poset and the cohomology groups. 
 
\subsection{The homotopy type of the complement}\label{sec3p1} First we look at arrangements in $S^1$. An arrangement in $S^1$ consists of $n$ copies of $S^0$, i.e. $2n$ points. The tangent bundle complement of such an arrangement is homeomorphic to the infinite cylinder with $2n$ punctures. Thus we have the following theorem.

\bt{thm0s1c4}
Let $\A$ be an arrangement of $0$-spheres in $S^1$. If $|\A| = n$ then 
\[M(\A) \simeq \bigvee_{2n+1} S^1.\] \et

From now on we assume that all our spheres are simply connected. We say that two arrangements are \textit{combinatorially isomorphic} if their corresponding face posets and intersection posets are isomorphic. \par

Let $\A$ be an arrangement of pseudo-circles in $S^2$. 
Then as a consequence of the \textit{Levi's enlargement lemma} \cite[Proposition 6.4.3]{ombook99} there exists a pseudo-circle $S_0\notin \A$ such that it is the equator with respect to the given antipodal map and it meets every member of $\A$ in exactly one point. 
Let $S_0^+, S_0^-$ denote the connected components of $S^2\setminus S_0$. 
As the equator $S_0$ intersects every $S\in \A$ generically, $S\cap S^{+}_0$ is a pseudo-line in $S^{+}_0 \cong \R^2$, respectively for $S^-_0$. 
Denote by $\A^+ := \A|S_0^+$, $\A^- := \A|S_0^-$ the pseudo-line arrangements in the respective hemispheres. 
Then $\A^+$ and $\A^-$ are combinatorially isomorphic arrangements of pseudo-lines.\par 

For the rest of the section we assume that $\A$ is centrally symmetric and there exists a hypersphere $S_0$ in general position such that the restriction of $\A$ to both the hemispheres (obtained by deleting $S_0$) result in combinatorially isomorphic, locally-flat pseudo-hyperplane arrangements (all the intersections are locally flat). 
This assumption motivates the following definition. 

\bd{defNew} A sphere arrangement $\A$ in $S^l$ is said to be \textit{mirrored} if it is centrally symmetric and there exists a pseudo-sphere $S_0\notin \A$ such that restriction of $\A$ to one of the connected components of $S^l\setminus S_0$ is results in a pseudo-hyperplane arrangement combinatorially isomorphic to the restriction of $\A$ to the other connected component. \ed 

Note that all sphere arrangements in dimensions $1, 2$ are mirrored. Not all sphere arrangements in higher dimension are mirrored since the enlargement lemma fails in general. We refer the reader to \cite[Proposition 10.4.5]{ombook99} for an example of pseudo-sphere arrangement in $S^3$ such that there is no equator in general position. However our assumption is not too restrictive, for example, arrangements corresponding to realizable oriented matroids are mirrored. If $\A$ is a mirrored sphere arrangement then we denote by $\A^+$ the restriction of $\A$ to one of the hemispheres and by $\A^-$ the restriction to the other. We note here that a lot of properties of hyperplane arrangements that we are interested in are also true for pseudo-hyperplane arrangements. For example, if all the intersections are locally flat then the construction of the associated Salvetti complex is same as described in Section \ref{sec:TheSalvettiComplex} and it has the homotopy type of the tangent bundle complement \cite[Theorem 3.3.7]{deshpande_thesis11}\par 

Here are two well known facts that we need. 

\bl{lem3s1c5} If $(Y, A)$ is a CW pair such that the inclusion $A\hookrightarrow Y$ is null homotopic then $Y / A\simeq Y\vee SA$, where $SA$ is the suspension of $A$.\el

\bpr See \cite[Chapter 0]{hatcher02}.\epr

\bl{lem4s1c5} Let $\mathcal{B}$ be an essential and non-central arrangement of pseudo-hyperplanes in $\R^l$ with locally flat intersections. Then the cell complex which is dual to the induced stratification is regular and homeomorphic to a closed ball of dimension $l$. \el

\bpr See \cite[Lemma 9]{sal1} and \cite[Proposition 9]{salvetti93}. \epr

Let $C$ be a chamber of $\B$ and for any other chamber $D$ let $\F(C, D)$ denote the set of all those faces $F$ such that $F\subseteq H$ if and only if $H\notin R(C, D)$. Clearly all the chambers of $\B$ are in $\F(C, D)$ and it is a disconnected set. In fact, if $l+1\leq |R(C, D)|\leq |\A|$ then $\F(C, D)$ is just the set of all chambers. Moreover, note that if $F\in \F(C, D)$ and $G$ some face such that $F\subseteq \ol{G}$ then $G\in \F(C, D)$. Let $Q$ denote the regular cell complex dual to the stratification induced by $\B$ and let $\F(C, D)^*$ denote the dual of $\F(C, D)$. It is straightforward to verify that each connected component of $\F(C, D)^*$ is contractible and deformation retracts onto a subset of $\partial Q$. \par 

For a chamber $C$ of $\mathcal{B}$ let $\iota_C$ denote the inclusion that takes a dual cell $F^*$ to the cell $\langle F^*, F^*\circ C^*\rangle$ of the associated Salvetti complex. This inclusion is same as the one introduced in Theorem \ref{thm2s3c3}. In fact, most of Theorem \ref{thm2s3c3} is true for pseudo-hyperplane arrangements (with appropriate modification in claim 3). For two distinct chambers $C, C'$ of $\mathcal{B}$ we define the following subset of the associated Salvetti complex
\[I(C, C') := \mathrm{Im}(\iota_C)\cap \mathrm{Im}(\iota_{C'}). \]

\bl{lem5s1c5}The subset $I(C, C')$ is non-empty and disconnected.\el

\bpr First observe that a cell $\langle G^*, D^*\rangle \in I(C, C')$ if and only if $D^* =G^*\circ C^* = G^*\circ (C')^*$. Which means that there is no hyperplane $H$ which contains $G$ and separates $C$ and $C'$. None of the cells obtained from a hyperplane $H\in R(C, C')$ can belong to $I(C, C')$ and hence, since all the vertices of the Salvetti complex are in $I(C, C')$, the set is disconnected. In fact, $\F(C, C')^*\cong I(C, C')$.\epr




\bt{thm1c4s1} Let $\A$ be a mirrored sphere arrangement in $S^l$. Let $\A^+$ and $\A^-$ be the pseudo-hyperplane arrangements in the two hemispheres and $\mathcal{C}(\mathcal{\A}^+)$ be the set of chambers of $\A^+$. Then the tangent bundle complement 
\[M(\A) \simeq Sal(\A^-)\vee \bigvee_{|\mathcal{C}(\mathcal{\A}^+)|}S^l. \] \et

\bpr The assumption that the arrangement is mirrored implies that $\A^+$ and $\A^-$ are combinatorially isomorphic, non-central, essential pseudo-hyperplane arrangements with locally flat intersections. Let $C\in \mathcal{C}(\A^+)$ and let $Q$ denote the dual cell complex $(S^+_0, \F^*(\A^+))$. Define the map $\iota_C^+\colon Q \hookrightarrow Sal(\A)$ as:
\[ F^* \mapsto  \langle F^*, F^*\circ C^*\rangle.\]

\textbf{Claim 1}: The image of the map $\iota_C^+$ in $Sal(\A)$ is homeomorphic to $Q$. \vskip2mm

Observe that $\iota^+_C$ is the restriction of the map $\iota_C$ defined in Theorem \ref{thm2s3c3}, which is an embedding. Hence $\iota_C^+$ maps $Q$ homeomorphically onto its image.\par

Thus $\iota_C^+$ is the characteristic map which attaches the boundary $\partial Q$ to the $(l-1)$-skeleton of $Sal(\A^-)$. For notational simplicity let $j_C$ denote the restriction of $\iota_C^+$ to $\partial Q$. \vskip2mm

\textbf{Claim 2}: The image of $j_C$ is the boundary of an $l$-cell in $Sal(\A^-)$. \vskip2mm

Consider the subcomplex of $Sal(\A^-)$ given by the cells $\{\langle F^*, F^*\circ C^*\rangle ~|~ F\in \F^*(\A^-) \}$. By Lemma \ref{lem3s1c5} above this subcomplex is homeomorphic to the closed $l$-ball. The boundary of this closed ball is precisely the image of $j_C$. \par

The characteristic map $\iota_C^+$ is the extension of $j_C$ to the cone over $\partial Q$ (which is $Q$). Hence $j_C$ is null homotopic.
Let $C'$ be any other chamber.
In view of Lemma \ref{lem5s1c5} the intersection set $I(C, C')\subseteq Sal(\A^+)$ retracts onto the boundary $\mathrm{Im}(j_C)\cup \mathrm{Im}(j_{C'})$.
Then it follows from Lemma \ref{lem3s1c5} that $Sal(\A^-) \cup \mathrm{Im}\iota^+_C(Q)$ has the homotopy type of $Sal(\A^-)\vee S^l$. 
Repeating these arguments for every chamber of $\A^+$ establishes the theorem. \epr

We state the following obvious corollary for the sake of completeness. 

\begin{center}
 \begin{figure}[ht]
 \includegraphics[scale=0.45, clip]{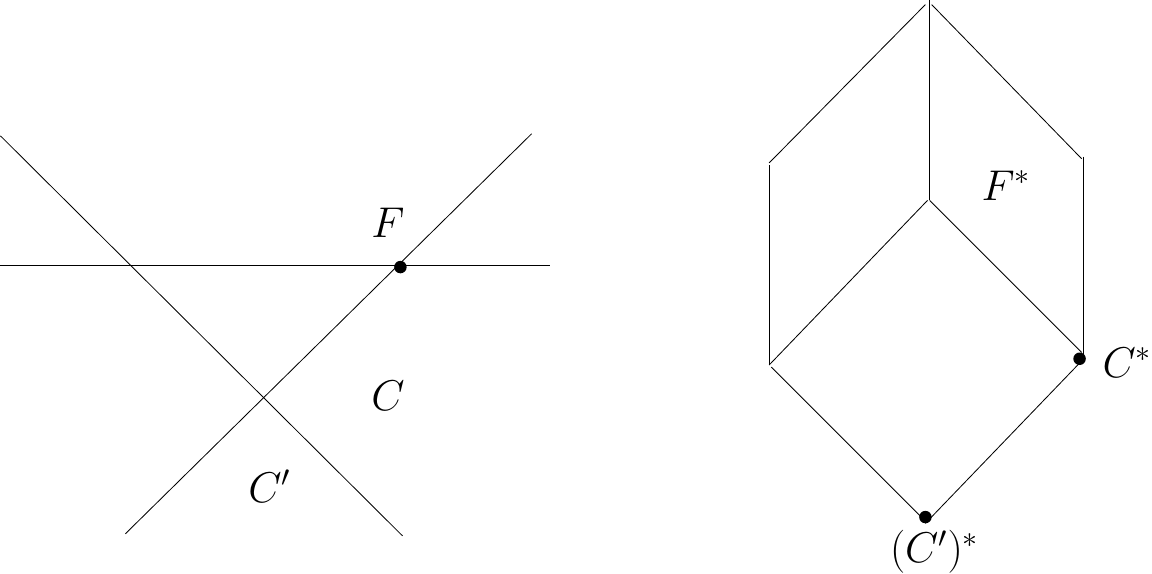}\caption{Restricted arrangement $\A^+$ and the associated dual complex $Q$.}\label{counterex}
  \end{figure}
\end{center}

\bc{cor1s1c4} Let $\A$ be a mirrored sphere arrangement in $S^l$. With notation as as before we have:
\[\pi_1(M(\A)) \cong \pi_1(M(\A^-)).  \] \ec

\be{ex1c4s1}
Consider the arrangement of $2$ circles in $S^2$ introduced in Example \ref{ex42}. It is clear that the arrangement $\A^-$ in this case is the arrangement of two lines in $\R^2$ that intersect in a single point. Hence 
\[M(\A)\simeq T^2 \vee S^2\vee S^2 \vee S^2 \vee S^2.\]
The Salvetti complex consists of four $0$-cells, eight $1$-cells and eight $2$-cells. The $2$-torus $T^2$ in the above formula corresponds to $M(\A^-)$. \ee

\be{ex3c4s1} Finally, consider the arrangement of three $S^2$s in $S^3$ that intersect like co-ordinate hyperplanes in $\R^3$. The $\A^-$ in this case is the arrangement of co-ordinate hyperplanes hence $Sal(\A^-) \simeq T^3$, the $3$-torus. This arrangement has $8$ chambers. So we have the following 
\[M(\A) \simeq T^3\vee \bigvee_8 S^3. \] \ee

\be{ex5c4s1} Consider the arrangement of three circles in $S^2$ that intersect in general position. This arrangement arises as the intersection of $S^2$ with the coordinate hyperplanes in $\R^3$. In this case $\A^-$ is the arrangement of three lines in general position. Figure \ref{counterex} shows $\A^+$ and the dual cell complex $Q$.
The intersection $I(C, C')$ (see Lemma \ref{lem5s1c5} above) contains the $2$-cell $\langle F^*, C^*\rangle = \langle F^*, (C')^*\rangle$. The boundaries of $\iota^+_{C}(Q)$ and $\iota^+_{C'}(Q)$ collapse to a point hence $Sal(\A^-)\cup \iota^+_{C}(Q)\cup \iota^+_{C'}(Q)$ has the homotopy type of $Sal(\A^-)\vee S^2\vee S^2$. In general the tangent bundle complement has the homotopy type of $Sal(\A^-)\vee \bigvee_7 S^2$.

\ee

\subsection{Cohomology of the Complement}\label{sec2p2}

We now establish a relationship between the cohomology of the tangent bundle complement and the intersection poset. Let $\A$ be a mirrored sphere arrangement in $S^l$, let $\A^+$ be the pseudo-hyperplane arrangement in the positive hemisphere. Let $L$ and $L^+$ denote the corresponding intersection posets. Observe that the map from $L$ to $L^+$ that sends $Y\in L$ to $Y|S^+_0 =: Y^+$ is one-to-one up to rank $l-1$. If $L_{l-1}$ and $L^+_{l-1}$ denote the sub-posets consisting of elements of rank less than or equal to $l-1$ then the previous map is a poset isomorphism. For notational simplicity we use $M^-$ for $M(\A^-)$. 

\bt{thm1nc4s1} With notation as above, we have the following 
\[ 
\ds \rank H^i(M, \Z) = \begin{cases}\sum_{\substack{Y\in L\\ \rank(Y) = i}} |\mu(S^l, Y)| &\text{for~} 0\leq i < l, \\ 
\sum_{Y\in L} |\mu(S^l, Y)| &\text{for~} i = l.\end{cases} 
\] 
Where $\mu$ is the M\"obius function of the intersection poset. \et

\bpr We use Theorem \ref{thm1c4s1} above and Theorem \ref{thm112} in order to prove the assertion by considering two cases. \vskip2mm
 \textbf{Case 1}: $1\leq i < l$. 
\begin{align*}
\rank(H^i(M)) &= \rank(H^i(M^-)) + \sum_{|\mathcal{C}(\A^+)|} \rank(H^i(S^l))\\
	&= \rank(H^i(M^-)) + 0 \\
	&= \sum_{\rank(Y^-) = i}(-1)^{\rank Y^-}\mu(\R^l, Y^-) \\
	&= \sum_{\rank(Y) = i} |\mu(S^l, Y)|.\end{align*}
The last equality follows from the fact that each $Y$ is a sphere of dimension $l-i$.\vskip2mm
\textbf{Case 2}: $i = l$. 
\begin{align*}
\rank(H^l(M)) &= \rank(H^l(M^-)) + \sum_{|\mathcal{C}(\A^+)|} \rank(H^l(S^l))\\
	&= \sum_{\rank(Y^-) = l} |\mu(Y^-)| + |\mathcal{C}(\A^+)|  \\
	&= \sum_{\rank(Y^-) = l} |\mu(Y^-)| + \sum_{Y^+\in L^+} |\mu(Y^+)| \\
	&= \sum_{Y\in L} |\mu(Y)|. \end{align*}
	
The third equality follows from the expression for the number of chambers of a hyperplane arrangement. The last equality is true because the number of rank $l$ elements in $L$ is twice the corresponding number in $L^-$. \epr

\br{rem4}
One might call the cohomology algebra $H^*(M(\A), \Z)$ the \emph{spherical Orlik-Solomon algebra}.
As stated in the introduction for mirrored arrangements the spherical OS-algebra is the direct sum of a graded and ungraded OS-algebras. 
The ungraded OS-algebra sits in the top dimension.
It should be interesting to figure out the structure of the spherical OS-algebra in the general case. 
Remark \ref{remint} implies that $H^*(M(\A), \Z)$ is a finitely generated $H^*(S^i, \Z)$ module for $1\leq i\leq l$.
\er
\section{The fundamental group}\label{sec4}

We turn our attention to the fundamental group of $M(\A)$. The aim of this section is to identify a class of sphere arrangements for which the word problem for $\pi_1(M(\A))$ is solvable. We do this by carefully analysing the so-called arrangement groupoid (in fact, the fundamental groupoid) and the path category of the associated Salvetti complex. This type of analysis goes back to the seminal work of Deligne in \cite{deli72}. 
\par 
In \cite{salvetti93} Salvetti introduced a certain class of cell complexes which were called \textit{metrical-hemisphere complexes} (MH-complexes for short). 
The combinatorial properties of the face poset of these MH-complexes resemble that of zonotopes.
He proves that these properties enables one to construct Salvetti complex similar to the description in the beginning of the previous section (see \cite[Section 2]{salvetti93}).
We should note here that MH-complexes are quite general and need not correspond to any type of arrangements.
However, Salvetti in his paper focuses on a certain class of MH-complexes which he calls MH*-complexes.
His main results are the following: Given an MH*-complex if the path category of the associated Salvetti complex admits a calculus of fractions then the fundamental group has solvable word problem \cite[Theorem 27]{salvetti93}. If in addition the dual of the given MH*-complex is a simplicial subdivision of a closed manifold then assuming an additional technical condition the universal cover of the Salvetti complex is contractible \cite[Theorem 33]{salvetti93}.
We make use of these ideas in the current context.
Instead of proving that the dual cell decomposition induced by a sphere arrangement is an MH*-complex we explicitly provide calculations that give a solution to the word problem.
As a result, most of the arguments in this section are straightforward generalizations of the results appeared in \cite{salvetti93}.
We have reproduced the proofs for the benefit of the reader.
\par 
An oriented $1$-skeleton of a regular CW-complex can be thought of as a directed graph without loops. The $0$-cells are the vertices and $1$-cells are the edges. A \textit{path} in such a cell complex is a sequence of consecutive edges and its \textit{length} is the number of edges. A \textit{minimal path} is a path of shortest length among all the paths with the same end points. Given a path we also define its \textit{initial vertex} and \textit{terminal vertex} in the obvious way. Finally, by a \textit{positive path} we mean a path all of whose edges are traversed according to their orientation.  We use the notation introduced in Section \ref{tbc}.

\bl{lem5s2c5} Let $\A$ be an arrangement of spheres in $S^l, l\geq 2$. Then any two minimal positive paths in the $1$-skeleton of $Sal(\A)$ that have same initial as well as terminal vertex are homotopic relative to $\{0, 1\}$. \el

\bpr Given two positive minimal paths $\alpha, \beta$ in $Sal(\A)$ with the initial vertex $\langle C^*, C^*\rangle$ and the terminal vertex $\langle D^*, D^*\rangle$ apply the retraction map $\psi$ defined in Theorem \ref{thm2s3c3}. 
Hence we get two paths in $\F^*$ and no two edges of $\alpha, \beta$ are sent to the same edge in $\F^*$. 
The conclusion follows from the observation that these image paths are contained in $\iota_C(\F^*)$, which is simply connected.\epr

Given an arrangement $\A$ we denote by $\G^+(\A)$ the associated \emph{positive category}, i.e., the category of directed paths in the Salvetti complex $Sal(\A)$. 
The objects of this category are the vertices of the Salvetti complex and morphisms are directed homotopy classes of positive paths (two such paths are connected by a sequence of substitutions of minimal positive paths). 
For a path $\alpha$ its equivalence class in $\G^+$ is denoted by $[\alpha]_+$.\par  

We denote by $\G(\A)$ the arrangement groupoid of $\A$. It is, in fact, the fundamental groupoid of the associated Salvetti complex. For a path $\alpha$ its equivalence class in $\G(\A)$ is denoted by $[\alpha]$. Since $\G(\A)$ is the category of fractions of $\G(\A)^+$ we denote by $J\colon \G(\A)^+\to \G(\A)$ the associated canonical functor. We refer the reader to \cite{gabzis67} for the relevant terminology from calculus of fractions. We drop the reference to $\A$ from $\G(\A)$ and $\G^+(\A)$ when the context is clear.\par 

For notational simplicity we write $F$ for $F^*$, i.e., we do not differentiate between a face in $\F(\A)$ and its dual cell in $\F^*(\A)$. Given a path $\alpha$ we denote by $(\pm a_1,\dots, \pm a_n)$ the sequence of $1$-cells in $Sal(\A)$ that are traversed by $\alpha$ either according to or opposite to their orientation depending on the sign.

\bt{thm1s1ch4}
Let $\A$ be a centrally symmetric arrangement of spheres in $S^l, l\geq 2$, then the associated canonical functor $J\colon \G^+\to \G$ is faithful on the class of minimal positive paths. \et

\bpr 
We already know that if $\alpha, \beta$ are two minimal positive paths with the same end points then $[\alpha] = [\beta]$ in $\G$. Hence it is enough to show that $[\alpha]_+ = [\beta]_+$. We argue on the lines of the proof of \cite[Theorem 20]{salvetti93}. Since each $S\in \A$ is centrally symmetric around the origin the antipodal map induces a fixed point free cellular action on the faces of $\A$.\par

Suppose $\alpha = (a_1,\dots, a_n)$ and $\beta = (b_1,\dots, b_n)$ are two minimal positive paths that start at $C$ and end at $D$. We proceed by induction on $n$, cases $n = 0, 1$ being trivial. Assume that the statement is true for all minimal positive paths with same end points and of length strictly less than $n$. If $a_1 = b_1$ then we are done by induction. Hence assume that $a := a_1$ and $ b := b_1$ are distinct and are dual to the hyperspheres $S_a, S_b$ respectively. \par

We have that $S_a, S_b\in R(C, D)$ (the set of hyperspheres separating $C$ and $D$) and that $S_a\cap S_b\cong S^{l-2}$. Recall that $S\in R(C, D)$ if and only if $C$ is contained in one of the connected components of $S^l\setminus S$ and $D$ is contained in the other. For $S\notin R(C, D)$ let $X_S(C, D)$ denote the closure of the connected component of $S^l\setminus S$ that contains both $C$ and $D$ and let $\h(C, D) = \bigcap_{S\notin R(C, D)}X_S(C, D)$. \par

\noindent\textbf{Claim 1}: The set $\h(C, D)$ is either connected or empty.\\
There are two cases: either $R(C, D) = \A$ or $R(C, D)\subsetneq \A$. In the first case it is clear that $\h(C, D)$ is empty since each $X_S(C, D)$ is empty. As for the second case; $\h(C, D)$ is an intersection of closed balls which is again a closed ball. See \cite[Lemma 5.1.9]{ombook99} and \cite[Page 217]{folk_lawr}. \par 

\noindent\textbf{Claim 2}: If $\h(C, D)\neq \emptyset$ then $\h(C, D)\cap S_a\cap S_b \neq \emptyset$. \\
Let $S_a^+, S_a^-$ (respectively $S_b^+, S_b^-$) denote the (closures of the) connected components of $S_a\setminus (S_a\cap S_b)$ (respectively $S_b\setminus (S_a\cap S_b)$). Without loss of generality assume that $S_a^+$ and $S_b^+$ intersect $\overline{C}$. This implies 
\[\h(C, D)\cap S_a^+ \neq \emptyset \neq \h(C, D)\cap S_b^+. \]
A similar argument using $\overline{D}$ establishes the claim. \par

Hence the set $\h(C, D)\cap S_a\cap S_b$ contains a codimension-$2$ face say $F^2$. Let $C'$ denote $F^2\ast C$. Let $\gamma_0$ be a minimal positive path from $C'$ to $D$. Also, there exist two minimal positive paths $\gamma_1, \gamma_2$ such that $\gamma_1$ starts at $a_1\circ C$ and $\gamma_2$ starts at $b_1\circ C$ such that both of them end at $C'$. Using this we can construct two new minimal positive paths $\eta = a_1\gamma_1\gamma_0$ and $\eta' = b_1\gamma_2\gamma_0$. The paths $\alpha, \eta$ are  minimal positive with the same end points and share the same first edge. Hence by induction, $[\alpha]_+ = [\eta]_+$ if $\gamma\neq 0$. For the same reasons $[\beta]_+ = [\eta']_+$. If $C'\neq D$ then the path $\gamma_0$ is of nonzero length and again by induction $[a_1\gamma_1]_+ = [b_1\gamma_2]_+$ implying $[\eta]_+ = [\eta']_+$. Now the transitivity of the equivalence relation proves the theorem when $\h(C, D)\neq\emptyset$. \par
The cases in which either $C' = D$ or $\h(C, D) = \emptyset$ can be treated similarly. \epr

For a centrally symmetric arrangement $\A$, we indicate by $[\mu(C\to D)]$ the unique equivalence class (in $\G^+(\A)$ or $\G(\A)$) determined by a minimal positive path from a chamber $C$ to another chamber $D$.\par

\bd{def3c3s5}
A sphere arrangement $\A$ in a $S^l$ is said to have \emph{the involution property} if there exists a graph automorphism $\phi\colon \F^*_1\to \F^*_1$ of the dual $1$-skeleton (considered as a graph) satisfying:
\begin{enumerate}
	\item $\phi$ is an involution (which induces involution on the vertices as well as the edges);
	\item for every vertex $C$, $d(C, \phi(C)) = \mathrm{max}\{ d(C, D) ~|~ D\in\Ch\}$;
	\item $d(C, \phi(C)) = d(C, D) + d(D, \phi(C))$ for every vertex $C$ and $D$. \end{enumerate}
\ed

\bl{lem2c4s1}
Let $\A$ be a centrally symmetric sphere arrangement in $S^l$ then $\A$ has the involution property. \el

\bpr The antipodal action on $S^l$ provides the required graph automorphism on the $1$-skeleton of the associated Salvetti complex. \epr 

The image of either a vertex or an edge under $\phi$ will be denoted by the superscript $\#$, for example, $C^\# := \phi(C)$.

\bl{lem1s5c3}
If $\A$ is a sphere arrangement with the involution property then
\begin{enumerate}
	\item $d(C, C^\#) = |\A|$ for all $C\in\Ch$;
	\item $d(C, D) = d(C^\#, D^\#)$ for all $C, D$.
\end{enumerate} \el

\bpr
Using the property (3) in Definition \ref{def3c3s5} we have:
\begin{align}
d(C, C^\#) &= d(C, D) + d(D, C^\#), \label{con1} \\
d(C, C^\#) &= d(C, D^\#) + d(D^\#, C^\#). \label{con2}
\end{align}
Adding equations (\ref{con1}) and (\ref{con2}) we get 
\[ 2 d(C, C^\#) = 2 d(D, D^\#).\]

Without loss of generality assume that $d(C, C^\#) = |\A|-1$. Hence there is a hypersphere $S\in \A$ such that $C$ and $C^\#$ are on the same side with respect to $S$. Choose a chamber $D$ such that $S\in R(C, D)$ and is adjacent to $C^\#$. Hence $d(C, D) > d(C, C^\#)$, which contradicts the Equation (\ref{con1}) above. Consequently no such $S$ exists. We call the number $d(C, C^\#) = |\A|$, the \emph{diameter of $S^l$} (with respect to $\A$). \par
Now subtracting $d(D, D^\#) = d(D, C^\#) + d(C^\#, D^\#)$ from (\ref{con1}) we get 
\[ 0 = d(C, D) - d(C^\#, D^\#)\]
which proves (2).\epr 

This involution also preserves the positive equivalence on paths as proved in the next lemma. 

\bl{lem2c3s5}
If $\A$ is a sphere arrangement  with the involution property then the involution $\phi$ induces a functor on $\G^+$ which is also an involution.\el 

\bpr
There is a bijection between the set of edge-paths of $\F_1^*$ and the set of all positive paths in $Sal(\A)_1$ given by $\langle F, C\rangle\mapsto F$. Extend the given involution to $Sal(\A)_1$ by sending $\langle F, C\rangle$ to $\langle F^\#, C^\#\rangle$. 
Under this involution a positive path $\alpha = (a_1,\dots,a_n)$ goes to a positive path 
\[\alpha^\# := (a_1^\#,\dots,a_n^\#).\]

If $\gamma_1, \gamma_2$ are two minimal positive paths contained entirely in the boundary of a $2$-cell of $Sal(\A)$ then so are $\gamma_1^\#, \gamma_2^\#$. Therefore $[\gamma_1]_+ = [\gamma_2]_+\Rightarrow [\gamma_1^\#]_+ = [\gamma_2^\#]_+$. \epr

We define a \emph{`positive'} loop based at $\langle C, C\rangle$ as 
\[ \delta(C) := \mu(C\to C^\#)\mu(C^\#\to C). \]
Note that the equivalence class of $\delta(C)$ in both $\G^+$ and $\G$ is unique.  
By $\delta^k(C)$ we mean that the positive loop is traversed $k$ times according to the orientation of the edges if $k>0$ and in the reverse direction if $k<0$. We will say that a positive path $\alpha$ begins (or ends) with a positive path $\alpha'$ if and only if $\alpha = \alpha'\beta (= \beta\alpha')$ for some positive path $\beta$.

\bl{lem3c3s5}
Let $\A$ be a sphere arrangement with the involution property and $\alpha$ be a positive path from $C$ to $D$. Then: 
\begin{enumerate}
	\item $[\alpha][\mu(D\to D^\#)] = [\mu(C\to C^\#)][\alpha^\#]$,
	\item if for a chamber $D'$, $\beta$ is some positive path from $C$ to $D'$ then $\alpha \delta^n(D)$ begins with $\beta$,
	\item if $[\gamma]\in \G(C, D)$ then there exists $n\in \N$ and a positive path $\gamma'$ such that 
	\[[\gamma] = [\delta^{-n}(C)][\gamma']. \] 
\end{enumerate} \el

\bpr
For ($1$) we use induction on the length of $\alpha$. In fact, it is enough to assume that $\alpha = \mu(C\to C_1)$ such that $d(C, C_1) = 1$. Thus:
\begin{eqnarray*}
\alpha \mu(C_1\to C_1^\#) &=& \mu(C\to C_1)\mu(C_1\to C_1^\#)  \\
			   &\stackrel{+}{\sim}& \mu(C\to C_1) \mu(C_1\to C^\#) \mu(C^\#\to C_1^\#) \\
			   &\stackrel{+}{\sim}& \mu(C\to C^\#) \mu(C^\#\to C_1^\#)  \\
			   &\stackrel{+}{\sim}& \mu(C\to C^\#) \alpha^\# .
\end{eqnarray*}
By the same arguments, the following stronger statement is true:

\begin{equation} [\alpha][\delta^k(D)] = [\delta^k(C)][ \alpha], \quad k\geq 1. \label{eq1c3s5} \end{equation}

For ($2$), let $\beta = (b_1,\dots, b_n)$ where $b_i$ is an edge from $B_{i-1}$ to $B_i$ ($B_0 = C, B_n = D'$). Observe that $\beta \mu(B_n\to B_{n-1}^\#) = (b_1,\dots, b_{n-1})\mu(B_{n-1}\to B_{n-1}^\#)$. 
By induction on $n$ assume that there exists a positive path $\eta$ from $B_{n-1}$ to $D$ such that - 
\[(b_1, \dots, b_{n-1})\eta = \alpha \delta^{n-1}(D). \]
Using ($1$), we get 

\[\beta \mu(B_n\to B_{n-1}^\#) \eta^\# = (b_1,\dots, b_{n-1})\eta \delta(D) = \alpha \delta^n(D) \]
which proves ($2$). \par
Let $\gamma$ be a path from $C$ to $D$. Assume $\gamma = (\epsilon_1 a_1,\dots, \epsilon_n a_n)$ where $\epsilon_i \in \{\pm 1\}$ for every $i$.\\ 
Let $A_i$ be the terminal vertex of $(\epsilon_1a_1,\dots,\epsilon_ia_i)$. 
Set $k = |\{1\leq i\leq n | \epsilon_i = -1  \}|$, we prove ($3$) by induction on $k$. The case $k = 0$ is clear since it means that $\gamma$ is a positive path. Assume that the statement is true for $k-1$. Now the general case; there exists an index $j$ such that $\epsilon_1 = \cdots = \epsilon_{j-1} = 1$ and $\epsilon_j = -1$. We have 
\begin{eqnarray*}
\delta(C)\gamma &=& \mu(C\to C^\#)\mu(C^\#\to C)(a_1,\dots,a_{j-1}, -a_j,\epsilon_{j+1}a_{j+1},\dots,\epsilon_n a_n)  \\
		&\stackrel{+}{\sim}& \mu(C\to C^\#)a_1^\#\mu (A_1^\#\to A_1) (a_2,\dots,a_{j-1},-a_j,\epsilon_{j+1}a_{j+1},\dots,\epsilon_n a_n)  \\
		&\stackrel{+}{\sim}& \mu(C\to C^\#)a_1^\#\cdots a_{j-1}^\#\mu(A_{j-1}^\#\to A_{j-1}) (-a_j,\epsilon_{j+1}a_{j+1},\dots,\epsilon_n a_n) \\
		&\stackrel{+}{\sim}& \mu(C\to C^\#)a_1^\#\cdots a_{j-1}^\# \mu(A_{j-1}^\#\to A_j) (\epsilon_{j+1}a_{j+1},\dots,\epsilon_n a_n)  \\
		&\stackrel{+}{\sim}& \delta^{1-n}(C) \gamma' ~(\hbox{by induction hypothesis}). 
\end{eqnarray*}

where $\gamma'$ is a positive path. Hence $[\gamma] = [\delta^{-n}(C)]\gamma'$.
\epr

Recall that \cite[Section 0.5.7]{stillwell80} the \emph{word problem} for a group $G$ is the problem of deciding whether or not an arbitrary word $w$ in $G$ is the identity of $G$. The word problem for $G$ is \emph{solvable} if there exists an algorithm to determine whether $w = 1_G$ or equivalently, if there exists an algorithm to determine when two arbitrary words represent the same element of $G$. 

\bt{thm3c3s5}
Let $\A$ be a sphere arrangement with the involution property. If the canonical functor $J\colon\G^+(\A)\to \G(\A)$ is faithful then the word problem for $\pi_1(M(\A))$ is solvable. \et

\bpr
Let $[\alpha], [\beta]$ be two loops in $\pi_1(Sal(\A))$ based at a vertex $\langle C, C\rangle$. Then according to Lemma \ref{lem3c3s5} there is a finite algorithm to write -
\begin{equation*}
[\beta] = [\delta^{-k}(C)] [\beta'],\quad [\alpha] = [\delta^{-k}(C)] [\alpha'] 
\end{equation*}
where $\beta', \alpha'$ are positive loops based at $\langle C, C\rangle$. Hence, $[\alpha] = [\beta]$ if and only if $[\alpha']_+ = [\beta']_+$. The theorem follows because there are only finitely many positive paths of given length to choose from. \epr

\br{rem2}
For simplicial arrangements of hyperplanes Deligne showed that the functor $J$ is injective \cite[Proposition 1.19]{deli72}. In particular it was used to prove that the word problem and the conjugacy problem for Artin groups are solvable \cite[Section 4.20]{deli72}. This was further generalized by Salvetti to include the case of simplicial arrangements of pseudohyperplanes in \cite[Theorem 31]{salvetti93}.
\er


\section{Arrangements of Projective Spaces}\label{sec5}
We now look at arrangements in projective spaces. Without loss of generality let $S^l$ be the unit sphere in $\R^{l+1}$ and $\p^l$ be the projective space with $a\colon S^l\to \p^l$ being the covering map. We consider a finite collection of subspaces that are homeomorphic to $\p^{l-1}$. We define the projective arrangements as follows. 

\bd{def3c4s1} A finite collection $\A = \{H_1,\dots, H_n\}$ of codimension-$1$ projective spaces is called an \emph{arrangement of projective spaces} (or a \emph{projective arrangement}) if and only if $\tilde{\A} := \{a^{-1}(H)~|~ H\in \A \}$ is a centrally symmetric arrangement of spheres in $S^l$. \ed

As a result of the projective version of the topological representation theorem, the following is an equivalent definition of projective arrangements. 
A finite collection of codimension-$1$ projective spaces (tamely embedded) in $\p^l$ is a projective arrangement if the non-empty intersection of every sub-collection is locally flat and homeomorphic to a lower-dimensional projective space. \par 

The homotopy type of the tangent bundle complement associated to a projective arrangement is easier to understand because of the antipodal action. 

\bt{thm2c4s1}
Let $\A$ be a projective arrangement in $\p^l$ and $\tilde{\A}$ be the corresponding centrally symmetric sphere arrangement in $S^l$. Then the antipodal map on the sphere extends to its tangent bundle and 
\[M(\A) \cong M(\tilde{\A}) / ((x, v)\sim a(x, v)). \] \et

\bpr
If $(x, v)$ is a point in the tangent bundle of $S^l$ extend the antipodal map in the obvious way, $a(x, v) = (-x, -v)$. We now prove that the space $M(\tilde{\A})$ is a covering space of $M(\A)$. This follows from the fact that $a\colon TS^l\to T\p^l$ is a covering map for every $l$. \par 
Note that the antipodal map is cellular on the faces of the arrangement. 
Consequently it induces a cellular map on $Sal(\tilde{\A})$ by sending a cell $\langle F, C\rangle$ to $\langle a(F), a(C)\rangle$ and we get a cell structure for $Sal(\A)$. 
Hence $\pi_1(M(\tilde{\A}))$ is an index 2 subgroup of $\pi_1(M(\A))$.\epr

\be{ex4c4s1}Here is a simple example. Consider the arrangement of pair of diametrically opposite points $\{p_1, p_2, q_1, q_2 \}$ in $S^1$. The arrangement breaks $S^1$ into $4$ chambers $A_1, A_2, B_1, B_2$; the chambers $A_1, B_1$ are diametrically opposite to $A_2, B_2$ respectively. The antipodal action is cellular and it identifies the $0$-cells $p_i$'s to $p$, $q_i$'s to $q$ and the $1$-cells $A_i$'s to $A$, $B_i$'s to $B$. Giving us the `projective' arrangement of two $\p^0$'s in $\p^1$. Figure \ref{figproj} shows the Salvetti complex of the sphere arrangement on the left and the Salvetti complex of `projective' arrangement on the right.
\begin{center}
\begin{figure}[ht]
\includegraphics[scale=0.5, clip]{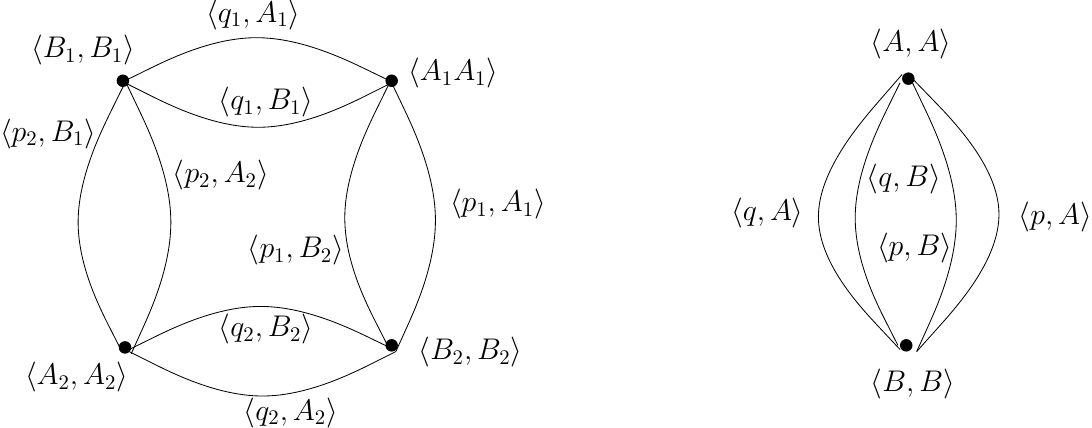}\caption{Spherical and projective Salvetti complexes}
\end{figure}\label{figproj}
\end{center}
\ee


Given a projective arrangement $\A$ let $J\colon \G^+\to \G$ denote the canonical functor between the positive category and the arrangement groupoid. For the corresponding (centrally symmetric) sphere arrangement $\tilde{\A}$ let $\tilde{J}\colon \tilde{\G}^+\to \tilde{\G}$ be the associated canonical functor. Recall that $\tilde{\A}$ has the involution property (Definition \ref{def3c3s5}) and that the antipodal action induces an `antipodal' functor on $\tilde{\G}^+$ (Lemma \ref{lem2c3s5}). Under this functor an object $C$ (which is a chamber) is mapped to $C^\#$ (its antipodal chamber) and a morphism $[\alpha]$ is mapped to $[\alpha^\#]$. 

\bl{lem5c4s1}
With the notation as above the following diagram commutes:
\[\begin{CD}
 \tilde{\G}^+ @>\tilde{J} >> \tilde{\G}\\
    @V\Phi^+ VV @VV \Phi V \\
    \G^+ @> J >> \G \\
\end{CD}\] 
where $\Phi^+$ identifies antipodal objects and morphisms and $\Phi$ is the covering functor for groupoids. \el

\bpr Follows from a simple diagram chase and the fact that $S^l$ is the universal cover of $\p^l$. \epr

An immediate consequence of the lemma is the following theorem. 

\bt{cor2c4s1} The restriction of $J$ to the class of minimal positive paths is faithful. Moreover if $J$ is faithful then the word problem for $\pi_1(M(\A))$ is solvable. \et

\bpr 
The first statement follows from the commutativity of the diagram in the previous lemma. If $[\alpha]_+$ is a class of minimal positive path in $\G^+$ then the class representing either of $\alpha$'s lift is also minimal positive in $\tilde{\G}^+$. If there are two distinct classes of minimal positive paths then first applying $\tilde{J}$ to their lifts in $\tilde{\G}^+$ and then applying  $\Phi$ results in producing two distinct classes of minimal positive paths in $\G$. By the same argument if $\tilde{J}$ is faithful then ${J}$ is also faithful. \par

Let $[\alpha]$ be a loop based at a vertex $C$ in $\G$. Let $[\tilde{\alpha}]$ be the class representing a lift of $\alpha$ which is a loop based at $\tilde{C}$ (a vertex in the fiber over $C$). Then by statement 3 in Lemma \ref{lem3c3s5} we have the following 
\[ [\tilde{\alpha}] = [\delta^{-n}(\tilde{C})][\tilde{\alpha}']\]
where $\delta^{-n}(\tilde{C}) = \mu(\tilde{C} \to \tilde{C}^\#)\mu(\tilde{C}^\#\to \tilde{C})$ and $\tilde{\alpha}'$ is a positive loop based at $\tilde{C}$. 
Since $\Phi$ is the covering functor, $\Phi([\delta^{-n}(\tilde{C})]) = [\delta^{-2n}(C)]$ here $\delta(C)$ is a positive loop based at $C$ which traverses every vertex twice. Let $[\alpha']$ be the image $\Phi([\tilde{\alpha}'])$, it represents a class of positive loop based at $C$. 
Note that choosing another lift of $\alpha$ based at the antipodal point $\tilde{C}^\#$ does not make any difference.
Hence we have proved that any loop in $Sal(\A)$ can be expressed as a composition of a `special loop' (which traverses each vertex a fixed number of times) and a positive loop. 
Now the same argument as in the proof of Theorem \ref{thm3c3s5} shows that the word problem for $\pi_1(M(\A))$ is solvable. 
\epr




\subsection*{Acknowledgments} This paper is based on a part of the author's doctoral thesis \cite{deshpande_thesis11}. The author would like to thank his supervisor Graham Denham for his support. The author would also like to acknowledge the support of the Mathematics department at Northeastern University for hosting a visit in 2011-12. The first version of this paper was written during that time. 


\begin{thebibliography}{99}

\bibitem{ombook99} A. Bj\"orner\ and\ M. Las Vergnas\ and\ B. Strumfels\ and\ N. White\ and\ G. Ziegler, 
{\it Oriented matroids}, second edition, Encyclopedia of Mathematics and its Applications, 46, Cambridge Univ. Press, 
Cambridge, 1999.

\bibitem{bz92} A. Bj\"orner\ and\ G. M. Ziegler, Combinatorial stratification of complex arrangements, J. Amer. Math. Soc. {\bf 5} (1992), no.~1, 105--149. 

\bibitem{bries73} E. Brieskorn, Sur les groupes de tresses [d'apr\`es V. I. Arnold], in {\it S\'eminaire Bourbaki, 24\`eme ann\'ee (1971/1972), Exp. No. 401}, 21--44. Lecture Notes in Math., 317, Springer, Berlin.

\bibitem{charney92} R. Charney, Artin groups of finite type are biautomatic, Math. Ann. {\bf 292} (1992), no.~4, 671--683. 

\bibitem{charney07} R. Charney, An introduction to right-angled Artin groups, Geom. Dedicata {\bf 125} (2007), 141--158.

\bibitem{dasdesh14} R. Das and\ P. Deshpande, Coxeter transformation groups and reflection arrangements in smooth manifolds. arXiv:1408.3921 [math.AT].

\bibitem{davisbook08} M. W. Davis, {\it The geometry and topology of Coxeter groups}, London Mathematical Society Monographs Series, 32, Princeton Univ. Press, Princeton, NJ, 2008.

\bibitem{deli72} P. Deligne, Les immeubles des groupes de tresses g\'en\'eralis\'es, Invent. Math. {\bf 17} (1972), 273--302.

\bibitem{deshpande_thesis11} P. Deshpande, {\it Arrangements of submanifolds and the tangent bundle complement}, PhD thesis, The University of Western Ontario, 2011. Electronic Thesis and Dissertation Repository. \url{http://ir.lib.uwo.ca/etd/154}.

\bibitem{folk_lawr} J. Folkman\ and\ J. Lawrence, Oriented matroids, J. Combin. Theory Ser. B {\bf 25} (1978), no.~2, 199--236. 

\bibitem{gabzis67} P. Gabriel\ and\ M. Zisman, {\it Calculus of fractions and homotopy theory}, Ergebnisse der Mathematik und ihrer Grenzgebiete, Band 35 Springer-Verlag New York, Inc., New York, 1967. 

\bibitem{gutkin86} E. Gutkin, Geometry and combinatorics of groups generated by reflections, Enseign. Math. (2) {\bf 32} (1986), no.~1-2, 95--110.

\bibitem{hatcher02} A. Hatcher, {\it Algebraic topology}, Cambridge Univ. Press, Cambridge, 2002. 

\bibitem{mandel82} A. Mandel, {\it Topology of oriented matroids}. Ph.D. Thesis, University of Waterloo (Canada). ProQuest LLC, Ann Arbor, MI, 1982.

\bibitem{miller87} D. A. Miller, Oriented matroids from smooth manifolds, J. Combin. Theory Ser. B {\bf 43} (1987), no.~2, 173--186.

\bibitem{orlik92} P. Orlik\ and\ H. Terao, {\it Arrangements of hyperplanes}, Grundlehren der Mathematischen Wissenschaften, 300, Springer, Berlin, 1992. 

\bibitem{rushing73} T. B. Rushing, {\it Topological embeddings}, Academic Press, New York, 1973. 

\bibitem{sal1} M. Salvetti, Topology of the complement of real hyperplanes in ${\bf C}\sp N$, Invent. Math. {\bf 88} (1987), no.~3, 603--618.

\bibitem{salvetti93} M. Salvetti, On the homotopy theory of complexes associated to metrical-hemisphere complexes, Discrete Math. {\bf 113} (1993), no.~1-3, 155--177. 

\bibitem{stillwell80} J. Stillwell, {\it Classical topology and combinatorial group theory}, second edition, Graduate Texts in Mathematics, 72, Springer, New York, 1993. 


\bibitem{dai1} D. Tamaki, Cellular Stratified Spaces I: Face Categories and Classifying Spaces. arXiv:1106.3772 [math.AT].

\end{thebibliography}
\end{document}